\documentclass[12pt,letterpaper]{amsart} 
\usepackage{palatino, euler, epic,eepic, amssymb, xypic, floatflt, microtype}
\usepackage{verbatim,color}

\usepackage{tikz-cd}

\input xy
\xyoption{all}

\newcommand{\point}{\vspace{3mm}\par \noindent \refstepcounter{subsection}{\thesubsection.} }
\newcommand{\tpoint}[1]{\vspace{3mm}\par \noindent \refstepcounter{subsection}{\thesubsection.} 
  {\bf #1. ---} }
\newcommand{\epoint}[1]{\vspace{3mm}\par \noindent \refstepcounter{subsection}{\thesubsection.} 
  {\em #1.} }
\newcommand{\bpoint}[1]{\vspace{3mm}\par \noindent \refstepcounter{subsection}{\thesubsection.} 
  {\bf #1.} }

\newcommand{\exercisedone}{ \vspace{2mm}}

\newcommand{\proj}{\mathbb P}

\newcommand{\C}{\mathbb{C}}

\newcommand{\Z}{\mathbb{Z}}

\newcommand{\R}{\mathbb{R}}

\newcommand{\cm}{{\mathscr{M}}}

\newcommand{\propernormal}{%
  \mathrel{\ooalign{$\lneq$\cr\raise.22ex\hbox{$\lhd$}\cr}}}

\newcommand{\defi}[1]{\underline{\smash{#1}}}

\usepackage{colonequals}
\usetikzlibrary{calc}
\usetikzlibrary{decorations.pathreplacing}
\usetikzlibrary{arrows,chains,matrix,positioning,scopes}

\makeatletter
\newcommand{\oset}[3][0ex]{%
  \mathrel{\mathop{#3}\limits^{
    \vbox to#1{\kern-2\ex@
    \hbox{$\scriptstyle#2$}\vss}}}}
\makeatother

\begin{document}

\pagestyle{plain}
\title{{\Large{The interpolation problem} } \\ \quad  \\ {\boldmath When can you pass a curve of a given type through $n$ random points in  space?}}
\author{Eric Larson}
\author{Ravi Vakil}
\author{Isabel Vogt}

\date{May 8, 2024.}
\begin{abstract}
The interpolation problem is a natural and fundamental question whose roots trace back to ancient Greece.  
The story is long and rich, with many chapters,  and a complete solution has been obtained only recently \cite{lv}.  Exploring it leads us  on a tour through a number of general themes in geometry.  This concrete problem motivates fundamental concepts such as moduli spaces and their properties, deformation theory, normal bundles, and more.  Questions about smooth  objects lead us to consider singular (non-smooth) objects, and in fact these smooth objects are studied by instead focusing on somehow simpler ``non-smooth'' objects, and then deforming them.
 \end{abstract}
\maketitle
\tableofcontents

{\parskip=12pt 

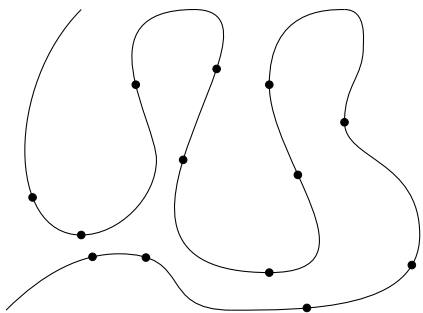
\begin{figure}[ht]
\begin{tikzpicture}
\draw (1, 4) .. controls (0, 3) and (0, 1) .. (1, 1);
\draw (1, 1) .. controls (1.5, 1) and (2, 1.5) .. (2, 2);
\draw (2, 2) .. controls (2, 2.5) and (1, 4) .. (2.5, 4);
\draw (2.5, 4) .. controls (4, 4) and (0.5, 0.5) .. (3.5, 0.5);
\draw (3.5, 0.5) .. controls (5, 0.5) and (3.5, 2) .. (3.5, 3);
\draw (3.5, 3) .. controls (3.5, 4) and (4.25, 4) .. (4.5, 4);
\draw (4.5, 4) .. controls (4.8, 4) and (4.75, 3.6) .. (4.75, 3.5);
\draw (4.75, 3.5) .. controls (4.75, 3.1) and (4.5, 3) .. (4.5, 2.5);
\draw (4.5, 2.5) .. controls (4.5, 2) and (5.5, 2) .. (5.5, 1);
\draw (5.5, 1) .. controls (5.5, 0) and (4, 0) .. (3, 0);
\draw (3, 0) .. controls (2, 0) and (2.5, 0.75) .. (1.5, 0.75);
\draw (1.5, 0.75) .. controls (1, 0.75) and (0.5, 0.5) .. (0, 0);
\filldraw (0.353, 1.5) circle[radius=0.05];
\filldraw (1, 1) circle[radius=0.05];
\filldraw (1.725, 3) circle[radius=0.05];
\filldraw (1.15, 0.71) circle[radius=0.05];
\filldraw (1.86, 0.7) circle[radius=0.05];
\filldraw (2.8, 3.21) circle[radius=0.05];
\filldraw (3.5, 0.5) circle[radius=0.05];
\filldraw (2.355, 2) circle[radius=0.05];
\filldraw (3.88, 1.8) circle[radius=0.05];
\filldraw (3.5, 3) circle[radius=0.05];
\filldraw (4.5, 2.5) circle[radius=0.05];
\filldraw (4, 0.03) circle[radius=0.05];
\filldraw (5.395, 0.6) circle[radius=0.05];
\end{tikzpicture}
  \caption{A curve passing through $13$ fixed randomly chosen points in $\R^2$}\label{f:fig1}
  \end{figure}

\section{A historical introduction to the interpolation problem}

Suppose that you are given $n$ distinct points in the plane $\R^2$, and you want to
pass a curve through them (see Figure~\ref{f:fig1}) --- not approximately as in a line of best fit
(which is an interesting problem in its own right), but exactly. Obviously,
given a pencil and enough hand-eye coordination you can do this.
As you increase the number of points, it's intuitive that any interpolating curve must also increase in ``complexity''.  To quantify this, we will flip the problem on its head and ask: if you fix some ``type'' of curve (which bounds its complexity) what is the maximum number of points that it interpolates?
Less precisely:  how ``bendy'' or ``flexible'' are curves of a given ``type''?
To motivate the precise statement of the problem, let's begin by considering some examples.

We could consider the ``type'' to be \emph{lines} (we will say that we are considering ``interpolation for lines'', implicitly specifying the ``type'' of curve we are considering) and ask: can one pass a {\em line} through the \(n\) points simultaneously?  If
$n >2$, then you'd have to be very lucky to be able to --- the points
would have to be collinear.  If $n=1$, it is very easy, and there are
infinitely many lines that will do the trick (and if $n=0$, it is
easier still).  Clearly $n=2$ is the important ``edge-case'' --- you can
manage it, but only one line does the job (Figure~\ref{f:fig2}).  This is the First Postulate in Euclid's elements \cite[Postulate 1, Book 1]{euclid}.

  \begin{figure}[ht]
\begin{tikzpicture}
\draw (0, 0) -- (4, 3);
\filldraw (1, 0.75) circle[radius=0.05];
\filldraw (3, 2.25) circle[radius=0.05];
\end{tikzpicture}
  \caption{There is a unique line through two
distinct    points in the plane (Euclid, c.\ 300 BCE)}\label{f:fig2}
  \end{figure}
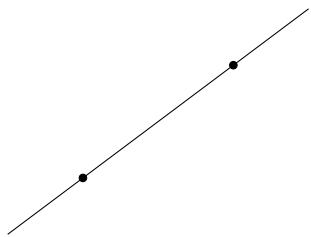

This is not the only solution to an interpolation problem that can be found in the {\em Elements}. It also contains the solution to the interpolation problem for circles.  By Proposition 5 of Book IV of \cite{euclid}, through \(3\) ``random'' points in \(\R^2\) there is a unique circle.  The fact that the circle is unique implies that \(3\) is the maximum number of random points, as then no circle passes through those three points and any fourth point off the circle.  (The four distinct points that {\em do} lie on a circle are the vertices of a cyclic quadrilateral, which is special because opposite angles add to \(\pi\)).

\begin{figure}[ht]
\begin{tikzpicture}[scale=.75]

\filldraw (0,0) circle[radius=0.06666666666666667];
\draw[densely dotted] (0,0) circle[radius=2];

\filldraw (1,1) circle[radius=0.06666666666666667];
\draw[densely dotted] (1,1) circle[radius=2];

\draw[densely dotted] (.5-2, .5+2) -- (.5, .5) -- (.5+2.5,.5 -2.5);

\filldraw(4, 1) circle[radius=0.06666666666666667];
\draw[densely dotted] (4,1) circle[radius=2];

\draw[densely dotted] (2.5, 3) -- (2.5, -2);

\draw (2.5, -1.5) circle[radius=2.915];

\end{tikzpicture}
\caption{The compass-and-straightedge construction of the unique circle interpolating \(3\) ``random'' points.}\label{f:compass_circle}
  \end{figure}
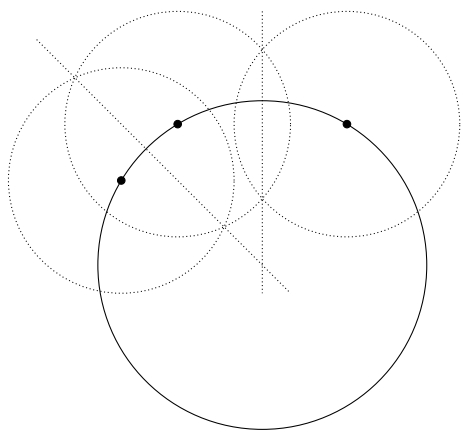

Notice that the word ``random'' is doing some serious work for us here.  If we were to choose three collinear points, then there would \textbf{not} exist a circle passing through these points.  And on the other hand, we can always find a circle interpolating as many {\em special} points as we want: simply cheat and choose the circle first and then pick points on it.  We use the term \defi{general} to describe ``random'' collections of points.  This is a notion coming from algebraic geometry that we will make more precise in the next section.  Roughly speaking, it means we want to avoid some unspecified special cases described by some algebraic conditions.  (See \S \ref{s:general} for more precision.)

\point ``The interpolation problem for a given type of curve'' is the following:  What is the {\em maximum} number of {\em general} points through which there always exists a curve of this type?

The answer to the interpolation problem can be viewed as a measure of the flexibility of this type of curve: each time you ask that the curve passes through an additional general point it must be flexible enough to accommodate.

We have already seen the answers to the interpolation problem for lines and circles: \(2\) and \(3\) respectively.

We make a brief but essential remark that so far our discussions are over the real numbers, but we will  soon be led to consider other fields.  So feel free to  continue to think of real numbers, but notice that often your ideas just automatically work in other fields (friendly ones such as $\mathbb{Q}$, or more dramatic ones such as in  \S \ref{rs}).  For example, the arguments above for lines and circles, suitably interpreted, will work over any field.

The following  classical examples will be important illustrations.

\bpoint{Polynomials taking specified values as an interpolation problem}
The  Lagrange interpolation problem asks: if I am looking for a polynomial $f(x)$ satisfying $n$ constraints
\[f(x_1)=y_1, \; f(x_2)=y_2, \; \dots, \; f(x_n) = y_n,\]
can I find such a polynomial of degree $d$?  We interpret this as an instance of our main problem, where the ``type'' is the graph of a degree \(d\) polynomial.  The ``generality'' condition will, in particular, ensure that \(x_i \neq x_j\) for \(i \neq j\).
When \(n = d+1\), we can actually write down an interpolating polynomial!  The key idea is that it is easy to find a polynomial satisfying \(f(x_1) = y_1\) and \(f(x_i) = 0\) for \(i \neq 1\):  by inspection,
\[f_1(x) = y_1 \frac{(x - x_2)(x-x_3) \cdots (x - x_{d+1})}{(x_1 - x_2)(x_1-x_3) \cdots (x_1 - x_{d+1})}\]
works.  We can  write down analogous \(f_2(x), \dots,  f_{d+1}(x)\), so  an interpolating polynomial is
\[f(x) = f_1(x) + f_2(x) + \cdots + f_{d+1}(x).\]
Since each of the polynomials \(f_i(x)\) has degree \(d\), the polynomial \(f(x)\) also has degree \(d\).
In fact, \(f(x)\) is the \emph{unique} interpolating polynomial of degree \(d\) --- if $g(x)$ were another then $f(x)-g(x)$ would be a polynomial of degree at most $d$ with $d+1$ roots $x_1$, \dots, $x_{d+1}$, so $f(x)-g(x) \equiv 0$.  (In terms of linear algebra, we have shown the nonvanishing of the Vandermonde determinant.)
 We have shown: 

\tpoint{The Lagrange Interpolation Theorem \cite{17, 27}}\label{t:lagrange} {\em The graph of a degree \(d\) polynomial interpolates \(d+1\) general points. }  \exercisedone

More explicitly: given \(d+1\) distinct  \(x_1, \dots, x_{d+1}\) and any \(y_1, \dots, y_{d+1}\), there exists a unique degree \(d\) polynomial \(f(x)\) such that 
\[f(x_1) = y_1, \dots, f(x_{d+1}) = y_{d+1}.\]
No such polynomial of degree \(d\) exists for general tuples \((x_1, y_1), \dots, (x_{d+2}, y_{d+2})\). 

The Lagrange interpolation theorem has been generalized by various other authors including \cite{4, 5, 15}.

\bpoint{Brief detour: Reed--Solomon error correcting codes}
\label{rs}
Theorem~\ref{t:lagrange} is a beautiful product of 18th century mathematics.  It is an equally beautiful 20th century insight of Reed and Solomon that this result has striking practical applications!
In their landmark 1960 paper \cite{reedsolomon}, Reed and Solomon showed  that the Lagrange Interpolation Theorem could be the basis for an efficient error correcting code.
Reed--Solomon error correcting codes are used to this day in most digital storage media, such as CDs or QR codes.

The basic problem is the following: suppose that Alice wants to send Bob \(n\) numbers \(p_1, \dots, p_n\).  But  Alice and Bob are communicating through a noisy channel that might transmit one of the numbers incorrectly.  
How can they build in redundancy to be able to detect, or even correct, errors in transmission?  
One approach would be to send each number twice or three times.  Then an error could be detected by noticing a mismatched pair, or corrected by taking the most common of the values.  While this  works, it is highly inefficient, since it requires sending \(2n\) or \(3n\) numbers. 
We can do substantially better using Lagrange Interpolation~\ref{t:lagrange}!

The key insight is that we can encode the numbers \(p_1, \dots, p_n\) as the coefficients of a polynomial \(f(x)\) of degree \(n-1\): 
\[f(x) = p_1x^{n-1} + p_2x^{n-2} + \cdots + p_{n-1}x + p_n.\]
By Theorem~\ref{t:lagrange}, the {\em values} of \(f(x)\) at \(n\) distinct inputs \(f(x_1), \dots, f(x_n)\) {\em determine} the polynomial \(f(x)\), and hence determine its coefficients \(p_1, \dots, p_n\) in a straightforward way.  Thus, instead of transmitting \(p_1, \dots, p_n\), Alice could instead transmit \(f(x_1), \dots, f(x_n)\).  So far, nothing has been gained: this is still a list of \(n\) numbers. 
(And we do have to ensure that nothing has been lost in this process!)
The added benefit  comes if we want to detect or correct errors.  Since \(n+1\) random pairs do not lie on the graph of a degree \(n-1\) polynomial, and errors can be expected to behave randomly, to detect a single error Alice simply sends {\em one} extra value of the polynomial \(f(x_{n+1})\).  If all goes according to plan, the \(n+1\) pairs of \(x_i\) and the value Bob
received for \(f(x_i)\) are special: they lie on the graph of a polynomial of degree \(n-1\).  Any \(n\) of them determine \(f\) and hence the desired \(p_1, \dots, p_n\).  
If an error occurs, Bob can recognize it: the transmitted pairs (almost certainly) do not lie on the graph of a polynomial of degree \(n-1\).  Since Bob won't know which is the erroneous pair to omit, in order to allow Bob to correct this error, Alice must send an additional value \(f(x_{n+2})\), so that Bob can recognize the correct polynomial as the one interpolating \(n+1\) pairs.
This is excellent news, since \(n+1\) or \(n+2\) numbers is {\em much} more efficient than \(2n\) or \(3n\).

\begin{figure}[h!]
    \centering
    \begin{minipage}{.45\textwidth}
\begin{center}
\begin{tikzpicture}
\draw [color=white, smooth, samples=100, domain=-1.3:1.5] plot(\x, {\x*\x*\x -
2.5*\x + (\x - 0.5) * (\x - 1) * (\x + 1)});
\draw [smooth, samples=100, domain=-1.5:1.7] plot(\x, {\x*\x*\x - 2.5*\x});
\filldraw (-1, 1.5) circle[radius=0.05];
\filldraw (0.5, -1.125) circle[radius=0.05];
\filldraw (1.0, -1.5) circle[radius=0.05];
\filldraw (-0.25, 0.609375) circle[radius=0.05];
\filldraw (-0.75, 1.45) circle[radius=0.05];
\filldraw (1.4, -0.756) circle[radius=0.05];
\end{tikzpicture}

{\small \(5\) points over-determine the graph of a degree \(3\) polynomial.}

\end{center}

    \end{minipage}
    \begin{minipage}{.48\textwidth}
        \begin{center}
\begin{tikzpicture}
\draw [smooth, samples=100, domain=-1.5:1.7] plot(\x, {\x*\x*\x - 2.5*\x});
\draw [densely dotted, smooth, samples=100, domain=-1.3:1.5] plot(\x, {\x*\x*\x -
2.5*\x + (\x - 0.5) * (\x - 1) * (\x + 1)});
\filldraw (-1, 1.5) circle[radius=0.05];
\filldraw (0.5, -1.125) circle[radius=0.05];
\filldraw (1.0, -1.5) circle[radius=0.05];
\filldraw (-0.25, 0.609375) circle[radius=0.05];
\filldraw[color=red] (-0.75, 2.0) circle[radius=0.05];
\filldraw (1.4, -0.756) circle[radius=0.05];
\end{tikzpicture}

{\small A single error can be corrected, since the other \(4\) points still over-determine a degree \(3\) polynomial.}

\end{center}
    \end{minipage}
    \caption{Illustration of the Reed--Solomon error correction protocol}
    \label{fig:rs}
\end{figure}
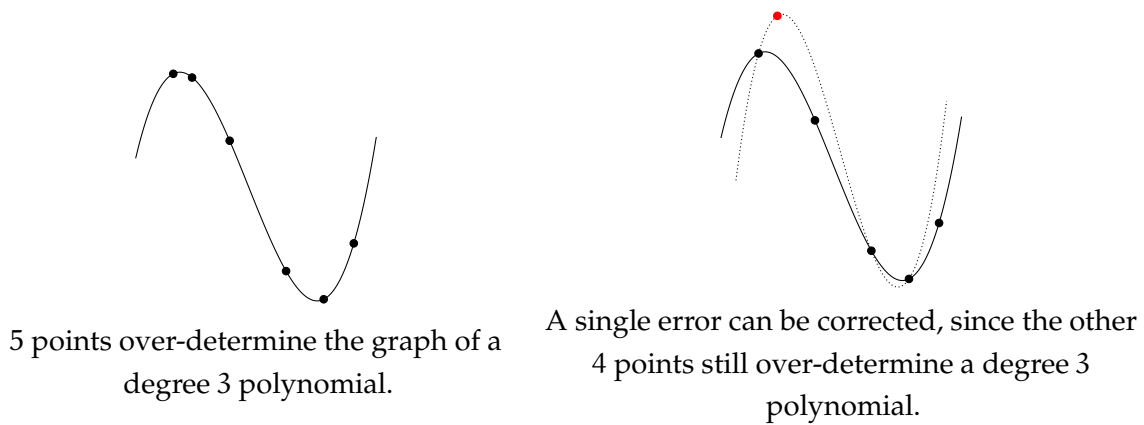

Of course, there are quite a few technical details to pay attention to here.  For example, if  we want this to actually be efficient, the value \(f(x_i)\) needs to be approximately the same number of bits as \(x_i\).  If the \(x_i\) are real numbers, this will not be the case!  This is seemingly fatal, but we get around the problem by instead working over a finite field \(\mathbb{F}_q\).  In this way, we are naturally led to consider the interpolation problem over more exotic fields than $\R$.  The notion of ``random'' or ``general'' is also more subtle in this case, which leads to some important practical caveats to the above rosy portrayal.

\bpoint{Conics}
Consider now \(5\) general points $(x_1, y_1)$, \dots, $(x_5, y_5)$ in $\R^2$.    Lagrange
interpolation ensures  that we can pass the graph of a degree $4$ polynomial
$y-f(x)=0$ through them.  But if our polynomial does not have to be a graph, we can manage something of lower degree
--- in fact, of degree $2$ (see Figure~\ref{f:conic}).  
We'll do this in two different ways here,
which will give us two different insights.  Later (\S \ref{s:conics2}) we will give yet another proof to highlight a technique employed in the proof of Main Theorem~\ref{t:main}.

\begin{figure}[ht]
\begin{tikzpicture}
\begin{scope}[rotate=30]
\draw (0, 0) ellipse (3 and 2);
\filldraw (-2.5, 1.1055415967851332) circle[radius=0.05];
\filldraw (-1.5, 1.7320508075688772) circle[radius=0.05];
\filldraw (-0.5, -1.9720265943665387) circle[radius=0.05];
\filldraw (2, -1.4907119849998598) circle[radius=0.05];
\filldraw (1, 1.8856180831641267) circle[radius=0.05];
\end{scope}
\end{tikzpicture}
  \caption{There is a (unique) conic through five generally chosen
    points (Pappus, c.\ 340 CE)}\label{f:conic}
  \end{figure}
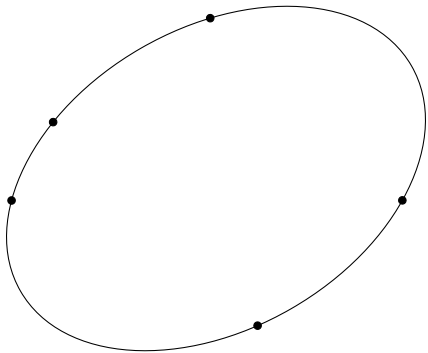

First, we can try to find a quadratic equation in two variables
that vanishes at these \(5\) points.  In other words, we are looking for
$6$ numbers  $a$, $b$, $c$, $d$, $e$, $f$ such that each of our \(5\) points $(x_i, y_i)$ 
lie on the curve described by the equation
\begin{equation}\label{eq:con}
a x^2 +  bxy + c y^2 + d x + e y + f =0.\end{equation}
The set of solutions \((x, y)\) to such an equation is called a \defi{conic},  so this is our chosen ``type'' of curve.
Since the conic is determined by these coefficients \(a, b, \dots, f\), we can try to find an interpolating conic by solving for the coefficients. 
In order to interpolate the 5 points that we start with, these coefficients must satisfy the \(5\) ``interpolation conditions'', one equation for each point \((x_i, y_i)\) to be on the conic:
\[a x_i^2 +  bx_iy_i + c y_i^2 + d x_i + e y_i + f = 0.\]
The key insight is that each is one {\em homogeneous linear} equation in the coefficients (since \(x_i^2\), \(x_iy_i\), etc are just scalars).
We can therefore rewrite the entire problem as
\[\begin{pmatrix} x_1^2 & x_1y_1 & y_1^2 & x_1 & y_1 & 1 \\
x_2^2 & x_2y_2 & y_2^2 & x_2 & y_2 & 1 \\
x_3^2 & x_3y_3 & y_3^2 & x_3 & y_3 & 1 \\
x_4^2 & x_4y_4 & y_4^2 & x_4 & y_4 & 1 \\
x_5^2 & x_5y_5 & y_5^2 & x_5 & y_5 & 1 \\
\end{pmatrix}  \begin{pmatrix} a \\ b \\ c \\ d \\ e \\ f\end{pmatrix} = \begin{pmatrix} 0 \\ 0 \\ 0 \\ 0 \\ 0 \\ 0\end{pmatrix}.\]
An interpolating conic has equation with coefficients \(a, b, c, d, e ,f\) that lie in the kernel of the above \(5 \times 6\) matrix.  Since a matrix with more columns than rows \textit{always} has a nontrivial kernel, we can always find a conic through \(5\)  points!  It may look like there are infinitely many solutions, but if the points are general, then the kernel will be \(1\)-dimensional, and so the possible \(a, b, c, d, e ,f\) will all be scalar multiples of each other.  In particular, they determine the same conic, which only depends on where the equation vanishes, and hence is insensitive to scaling the entire equation by a nonzero scalar.  We cannot hope to find a conic through \(6\) general points, since a \(6\times 6\) matrix of this shape will (in general) have only a trivial kernel. 

It will be helpful to step back and reinterpret this result.
How could we predict that the number of general points was \(5\)?  Each interpolation condition gives rise to a single equation.  On the other hand the number of unknowns \(a, b, \dots, f\) we are solving for is \(6\), not \(5\).
It looks like we have one too many unknowns.  But multiplying the
equation of a conic by a nonzero scalar yields the same curve.
In other words, plane conics are parameterized not by $\R^6$, but by (nonzero) elements of \(\R^6\) up to common scaling by a nonzero scalar.
In other words, the parameter space of conics is \defi{real projective $5$-space} \(\R \proj^5\), where
\[\R \proj^n \colonequals \frac{\R^{n+1} \smallsetminus \{0\}}{\R^\times}.\]
We can similarly define \(\C\proj^n\), etc., and when we wish to remain ambiguous about the field, we will write \(\proj^n\).
(Notice how we are forced to consider the notion of projective
geometry, just as the ancients were!)
Hence we observe that \(5\) --- the answer to the interpolation problem for plane conics --- is nothing other than the dimension of the parameter space of plane conic curves.

\epoint{Plane conics, take two} \label{s:parametrically}Now let us instead think of the problem parametrically.
The right way of posing the question turns out to be the following:  Are there three degree two polynomials $X(t) = a t^2+bt + c$, $Y(t) =
d t^2 + e t +f$, and $Z(t) = g t^2 + h t + i$ such that the
parameterized curve
\begin{equation}\label{eq:parcur}(X(t)/Z(t),
Y(t)/Z(t))\end{equation} passes through the \(5\) given points?
 This turns out to be the same question, because such parameterized
 curves are precisely the plane conics (although we will not justify that here).  It looks like we have $9$
 unknowns $a$, \dots, $i$, which is well over our \(5\) constraints.
 We realize quickly that we can knock this down to $8$ unknowns,
 because
 we can multiply \(a, \dots, i\) by a common nonzero scalar, but we are still off by three unknowns.  It turns out that
 there is a 3-dimensional way of parameterizing a conic, and that is
 where these 3 extra choices come from.  (This $3$ is precisely the
 dimension of the automorphism group of $\proj^1$.  This was also
 well-known to the classical Greeks, although not in this language.)

\bpoint{Higher degree plane curves}
\label{s:higherdegree}A \defi{degree \(d\) plane curve} is the set of solutions \((x, y)\) to an equation of degree \(d\) in \(x\) and \(y\).  Since the curve only determines the coefficients up to common scaling, the parameter space of degree \(d\) plane curves is \(\R\proj^N\), where \(N+1\) is the number of  monomials of degree \(d\) in \(x\) and \(y\).  Determining \(N\) is a  combinatorial exercise: \(N + 1 = \frac{(d+1)(d+2)}{2}\), and so \(N = \frac{d(d+3)}{2}\).   (For example, the equation of a cubic has $10$ coefficients, see \eqref{eq:cubcur}, so plane cubics are parameterized by $\R \proj^9$.)  At this point you probably agree that linear-algebraic considerations imply that the number of general points that a degree \(d\) plane curve interpolates is \(\frac{d(d+3)}{2}\).  This is a classical theorem from 1750 CE of Cramer (of linear algebraic fame) \cite{8}.

There is more to extract from this example, and we can already see this in degree $3$.
There is a single cubic 
 \begin{equation}\label{eq:cubcur}
 a x^3 + b x^2 y + c x y^2 + d y^3 + e x ^2 + f xy + g y^2 + h x + iy + j = 0\end{equation}
passing through $9$ generally chosen points.
However, a similar calculation as in the case of parameterized conics (\S \ref{s:parametrically}) shows that there is \emph{not} a degree $3$ {\em parameterized curve} passing
through $9$ generally chosen points --- you can only hope to pass one through
$8$ generally chosen points.  
What is going on?

The answer is helpfully explained by 
passing from a real curve to the set of solutions to the same system of equations over the complex numbers.  Since \(\mathbb{C}\) has two real dimensions, this ``curve'' now looks like a real surface, known as a \defi{Riemann surface}.
For example, from \(\mathbb{R}\proj^1 = \mathbb{R} \cup \{\infty\}\), we pass to \(\mathbb{C} \proj^1 = \mathbb{C} \cup \{\infty\}\), which is topologically a sphere --- known as the \defi{Riemann sphere}.
Every
parameterized complex curve \eqref{eq:parcur} of any degree is parameterized by \(t \in \mathbb{C}\proj^1\), and hence is basically
the Riemann  sphere (except we remove finitely many points where
the denominators $Z(t)$ are $0$).
But it turns out that a cubic curve \eqref{eq:cubcur} is torus (a ``genus $1$ Riemann
surface''), again minus finitely many points.    
The key additional information here is
that our curves now have a notion of \defi{genus} (the number of holes),
and cubic curves can have genus  $0$ or genus $1$.  Genus $0$ cubic curves
can interpolate $8$ general points in the plane, and genus $1$ cubic curves
can interpolate $9$ general points in the plane.  (See
Figure~\ref{f:genus2} for a sketch of a genus $2$ Riemann surface.)

\begin{figure}[ht]
\begin{tikzpicture}
\draw (0, 1) .. controls (0, 0) and (2, 0) .. (3, 0);
\draw (3, 0) .. controls (6, 0) and (7.5, 1) .. (7.5, 2);
\draw (0, 1) .. controls (0, 1.5) and (1, 2) .. (1.5, 2);
\draw (1.5, 2) .. controls (2.25, 2) and (2, 1.75) .. (3, 1.75);
\draw (3, 1.75) .. controls (4, 1.75) and (4, 3) .. (5.5, 3);
\draw (5.5, 3) .. controls (6, 3) and (7.5, 3) .. (7.5, 2);
\draw (1.5, 1) .. controls (2, 0.5) and (2.5, 0.5) .. (3, 1);
\draw (2, 0.666666) .. controls (2.16666, 0.833333) and (2.33333, 0.833333) .. (2.5, 0.66666);
\draw (4.5, 1.5) .. controls (5, 1) and (5.5, 1) .. (6, 1.5);
\draw (5, 1.166666) .. controls (5.16666, 1.333333) and (5.33333, 1.333333) .. (5.5, 1.16666);
\end{tikzpicture}
  \caption{A genus $2$ Riemann surface}\label{f:genus2}
  \end{figure}
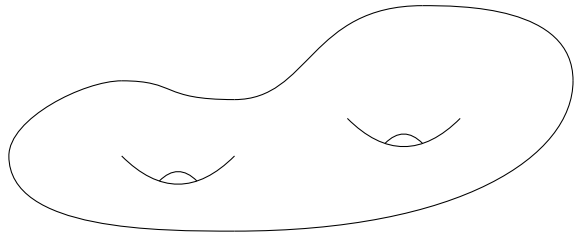

This allows us to refine our notion of ``type'' of plane curve: there is the degree (of the defining equation) and also the genus (of the Riemann surface over the complex numbers).
In fact the notion of genus can also be defined in a purely algebraic manner (not just in terms of ``holes of a Riemann surface''), and thus makes sense over any field.

\bpoint{All fields are now algebraically closed}
 For simplicity, from now on we will assume that {\em we are working over an algebraically closed field.}  (Until now this has not been an issue, since all of our thinking has been in terms of linear algebra, which works over any field.)

\bpoint{Interpolation for degree $d$ genus $g$ plane curves}
We can now state the full story of interpolation for degree $d$ curves in the plane.
Such curves can have genus $g$ anywhere between $0$ and $(d-1)(d-2)/2$
inclusive.
(A generally chosen  degree $d$ polynomial in two variables has genus
$(d-1)(d-2)/2$.
At the other extreme, a ``parameterizable'' degree $d$ curve, also known
as a \defi{rational curve}, has genus $0$.)

\tpoint{Theorem (interpolation for degree $d$ genus $g$ plane curves)} \label{t:intplane}
{\em Plane curves of degree $d$ and genus $g$ can
interpolate $3d+g-1$ generally chosen points in the plane.} \exercisedone

You can check that if \(g\) is equal to the largest possible value \((d-1)(d-2)/2\), then this is Cramer's theorem (\S \ref{s:higherdegree}).  If $g$ is equal to the smallest value (zero), you should be able to verify the result by following the parametrical proof by generalizing the argument for conics given in \S \ref{s:parametrically}.

\bpoint{What's next}\label{s:quadric_surface}
Our tour so far has taken us  through a vast territory and
chronology of mathematics, from millennia past through the late
nineteenth century.
This yielded a rather complete picture of the interpolation problem for curves in the plane.  

Where do we go from here?
One option might be to consider the solution sets of a single
equation in more variables.  This translates into studying the interpolation problem for
higher dimensional algebraic \defi{hypersurfaces}.
In this setting, the condition of interpolating each point is {\em still} a single homogeneous linear equation in the coefficients of the polynomial.
At this point, you probably feel confident that the number of general points such a hypersurface can interpolate is equal to the number of coefficients minus one (and know how to use linear algebra to prove it!).  To separate out a specific example that will be important later, consider \defi{quadric surfaces} in \(3\)-dimensional space; these surfaces correspond to the solutions \((x, y, z)\) to a single quadratic equation
\[
ax^2 + b xy + c xz + d y^2 + eyz + fz^2 + gx + hy + iz + j = 0.
\]
Since there are \(10\) coefficients, we conclude that quadric surfaces interpolate \(9\) general points.

While this will be an important calculation later on, it is not particularly novel --- once we understand plane curves, we can easily generalize to hypersurfaces of any dimension.
Instead of keeping the property of being described by the vanishing of a single equation (which led to interpolation for hypersurfaces), we will keep the property of being \(1\)-dimensional.  In other words, we will continue to study the interpolation problem for curves, but instead of curves in the plane, we consider curves in higher-dimensional space.  More precisely, we rethink of our curves in the plane as curves in the projective plane $\proj^2$  (so $x^2+x + y^2=1$ becomes $x^2+xz + y^2=z^2$ in ``projective coordinates''), and in our generalization, we consider curves in higher-dimensional projective spaces $\proj^r$. 

This will bring the problem forward into the twentieth century,
where the interpolation problem for curves has remained an active area of research
(including \cite{p6, lv, p5, p2, p1, p3, p4}).
More sophisticated algebraic geometry comes into play in order to understand what types of curves we will  have to consider.

  \section{More sophisticated algebraic geometry emerges: Brill--Noether theory}

 We now extend our discussion to consider maps of curves to  higher-dimensional projective spaces $\proj^r$.

  The notion of ``type'' of curve to consider will take on a new
  complication, as we will soon see.  We still focus on curves of a
  given genus $g$ (an intrinsic discrete invariant of curves). 
  
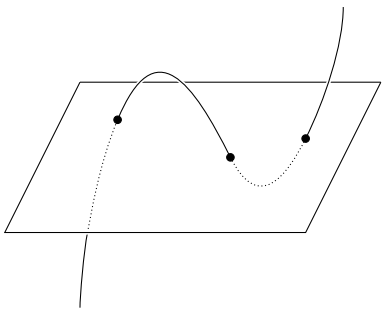
\begin{figure}[h!]
\begin{tikzpicture}
\draw (1, 3) -- (5, 3);
\draw[white, ultra thick] (1.5, 2.5) .. controls (2, 3.7) and (2.5, 3) .. (3, 2);
\draw[white, ultra thick] (4, 2.25) .. controls (4.25, 2.75) and (4.5, 3.5) .. (4.5, 4);
\draw (1, 0) .. controls (1, 0.15) and (1.05, 0.75) .. (1.1, 1);
\draw[densely dotted] (1.1, 1) .. controls (1.15, 1.25) and (1.25, 1.9) .. (1.5, 2.5);
\filldraw (1.5, 2.5) circle[radius=0.05];
\draw (1.5, 2.5) .. controls (2, 3.7) and (2.5, 3) .. (3, 2);
\filldraw (3, 2) circle[radius=0.05];
\draw[densely dotted] (3, 2) .. controls (3.375, 1.25) and (3.75, 1.75) .. (4, 2.25);
\filldraw (4, 2.25) circle[radius=0.05];
\draw (4, 2.25) .. controls (4.25, 2.75) and (4.5, 3.5) .. (4.5, 4);
\draw[white, ultra thick] (4, 1) -- (0, 1);
\draw (5, 3) -- (4, 1) -- (0, 1) -- (1, 3);
\end{tikzpicture}
    \caption{The degree of a curve in \(\proj^3\) is the number of times it meets a plane.}
    \label{f:degree}
\end{figure}

  Curves in $\proj^r$ also have a notion of degree, although it cannot
  quite be understood as the ``degree of defining equations''.
  Instead, it can be described as follows.  For simplicity, we work
  over the complex numbers.  The degree of a complex curve in
  $\proj^r$ is the number of points in which it meets a general hyperplane --- i.e., the locus defined by a single linear equation (see Figure~\ref{f:degree}).  Equivalently, its homology class is $d$ times
  the class of a line in $H_2(\C \proj^r, \Z)$.  (The class of a line is in
  fact a generator of $H_2(\C \proj^r, \Z) = \Z$.)  The equivalence of
  these  notions of degree is a special case of ``B\'ezout's Theorem''.

  So as in Theorem~\ref{t:intplane}, we have three invariants of our curves:  $r$ (corresponding to the ambient space $\proj^r$), $g$ (the genus), and $d$ (the degree).    In some sense, $g$ and $d$ measure the ``complexity'' of the curves in question.

  But this information is not enough to yield a good question, and
  therein lies an interesting story.  Curves with invariants $(g,r,d)$
  don't form just a set --- they form a {\em family}, or more
  precisely, a {\em moduli space}, which we call (for the purpose of this article only) $\cm(g,r,d)$.  We are led to the insight that to
  understand geometric objects of a given type, we should consider the
  {\em parameter space} or 
  {\em moduli space} of objects of that type.  For example,
  we have already seen that the family of plane conics \(\cm(0, 2, 2)\) is  a dense open subset of \(\proj^5\).  More generally, 
  curves
  in $\proj^2$ with given $g$ and $d$ form a family of dimension
  $3d+g-1$.  This requires making sense of the dimension of a family, but hopefully Theorem~\ref{t:intplane} gives some idea of how this must work.  (If we are working over $\C$, then this dimension is ``complex dimension'', for example.)  Moreover, this family has only one ``piece'' (``irreducible component'', see Figure~\ref{f:threepiece} for a ``definition'').

\begin{figure}[ht]
\begin{tikzpicture}
\draw (4.5, 1) -- (3, 4);
\draw[ultra thick, white] (3, 2) -- (6, 4);
\draw (3, 2) -- (6, 4);
\draw (0, 0) -- (1.875, 1.25);
\draw[densely dotted] (1.875, 1.25) -- (3, 2);
\draw (6, 3) .. controls (5.75, 2.75) and (5.7, 2.3) .. (6, 2);
\draw (6, 2) .. controls (6.7, 1.3) and (6.5, 0.5) .. (6, 0);
\draw[ultra thick, white] (1, 3) -- (2.5, 0);
\draw (1, 3) -- (2.5, 0);
\draw (2.5, 0) -- (4.5, 1);
\draw (3, 4) -- (1, 3);
\end{tikzpicture}
\caption{The ``irreducible components'' of a space are its ``pieces''.  For example, this space has three irreducible components.  Two have dimension one, and one has dimension two.}
\label{f:threepiece}
\end{figure}
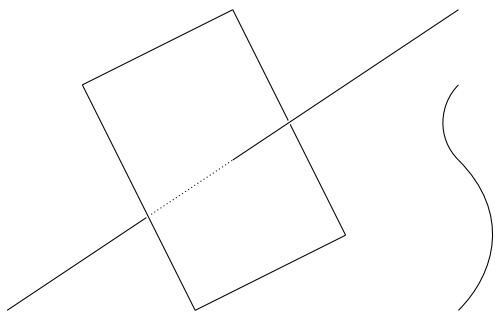

\epoint{Finally: The meaning of ``general''}
\label{s:general}We are finally ready to make precise our intuition of a ``general'' object of a certain type.   
If the objects naturally form some sort of algebraic moduli space (whose points correspond to the objects in question), then the phrase \defi{general object} means an object corresponding to a point in a dense open subset of the moduli space, whose complement is described by a bunch of algebraic equalities.  For example, a general quadratic in one variable $a x^2 + bx + c=0$ over the complex numbers has distinct roots, because the quadratics having a double root are those satisfying $b^2=4ac$.   The nice thing about this word is that we don't actually have to know what the algebraic equations describing the complement are --- we just have to know that such algebraic equations exist.    So this notion of ``general'' can be applied, for example, to any  type of object with an irreducible moduli space.  In this way, we can talk about, for example, a general curve of genus $g$, or a general $n$-tuple of points in $\proj^r$.

  \epoint{Dimensions of moduli spaces suggest how many points the corresponding  curves can interpolate} \label{s:countingpointconditions}The fact that the dimension of the moduli space $\cm(g,2,d)$ is $3d+g-1$ relates to Theorem~\ref{t:intplane}, and is much of the proof, as follows. (\defi{Base case  = Step $0$}) there is a $(3d+g-1)$-dimensional family of degree \(d\) and genus $g$ curves in $\proj^2$.  (\defi{Step $1$}) There is a $(3d+g-2)$-dimensional family of such curves through a fixed (general) point $p_1\in \proj^2$.
  (\defi{Step $2$}) There is a $(3d+g-3)$-dimensional family of such curves through two fixed (general) points $p_1, p_2$.   The pattern continues until:
  (\defi{Step $3d+g-1$}) There is a $0$-dimensional family of such curves through $3d+g-1$ fixed (general) points $p_1, p_2, \dots, p_{3d+g-1}$.  This last family is a finite (and, in fact, nonzero) number of points!     And finally: (\defi{Step $3d+g$}) there are no such curves through $3d+g$ fixed (general) points.    This suggests
  that such curves interpolate  precisely $3d+g-1$ general points, which is  correct.
  
\epoint{Moduli spaces of curves  in higher-dimensional projective spaces behave pathologically}
  In $\proj^r$, for $r>2$, aside from luckily chosen $g$ and $d$,
  the space $\cm(g,r,d)$  of curves with these invariants will have many different
  pieces (irreducible components) of many different dimensions (see again Figure~\ref{f:threepiece}).  Examples of this were first
  described by Mumford in his dramatic paper \cite{mumford}.  More generally, Mumford predicted that these spaces for
  various $g, r, d$ should satisfy ``Murphy's Law'', that they behave ``as badly as possible''.  In particular, 
  they can have any number of irreducible components, whose dimensions can differ
  by any given amount (and in fact can have arbitrary singularities in
  a well-defined sense), \cite{murphy}.  There is no reasonable hope of enumerating these families, let alone figuring out their dimensions.

  Here is one example, to give a sense of what can go wrong.  We will give two fundamentally  different families of degree $98$, genus $78$ curves in $\proj^{23}$.

\epoint{The first family}\label{s:first_family}
  The first family of curves is quite concrete (but perhaps not the first family you might think of).  Choose a degree $14$ homogeneous polynomial in three variables; this gives a curve in $\proj^2$ of degree $14$, and genus $\binom {13} 2 = 78$ by our discussion in \S \ref{s:higherdegree}.  Then map this $\proj^2$ in turn to $\proj^{23}$ by taking $24$ randomly chosen degree $7$ polynomials in $3$ variables.    There is enough information in our discussion so far to see that this forms a $974$-dimensional family.  (The idea:  There is a $(\binom {16} 2$-1)-dimensional family of such curves in $\proj^2$, and  $\binom 9 2$-dimensional choices for each of the degree $7$ polynomials, except multiplying them all by the same nonzero scalar gives the same map to $\proj^{23}$ so we subtract one; and we also need to mod out by automorphisms of $\proj^2$, so we subtract $8$ more.)   See Figure~\ref{f:crazycurve} for a sketch of the construction.  This family is maximal, in that there is no higher-dimensional irreducible family containing it. With a little thought, you may be able to interpret this as saying that if you take a  randomly chosen curve in this family, and try to deform it a little (wiggle it), it is still  a curve in this family (so in particular, it must continue to sit in one of these nonlinearly embedded projective planes in $\proj^{23}$).

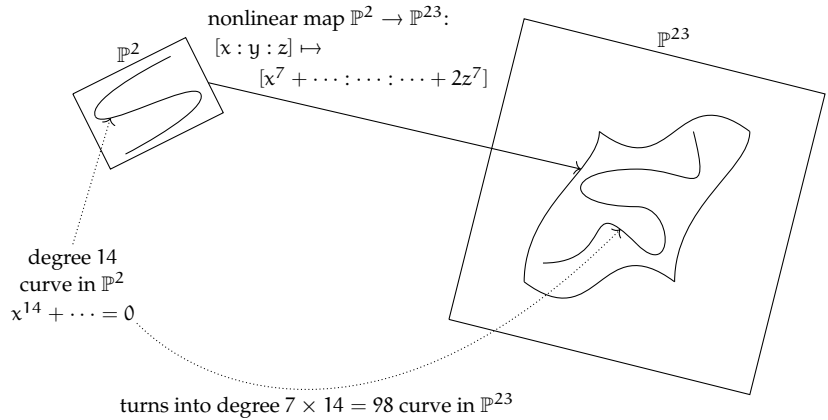
\begin{figure}[ht]
\begin{tikzpicture}
\draw (-1, 4) -- (0.5, 4.75) -- (1, 3.75) -- (-0.5, 3) -- (-1, 4);
\draw (4, 1) -- (8, 0) -- (9, 4) -- (5, 5) -- (4, 1);
\draw (0.3, 4.5) .. controls (-0.1, 4.3) and (-0.8, 3.9) .. (-0.7, 3.7);
\draw (-0.7, 3.7) .. controls (-0.6, 3.5) and (0.6, 4.2) .. (0.7, 4.0);
\draw (0.7, 4.0) .. controls (0.8, 3.8) and (-0.1, 3.3) .. (-0.3, 3.2);
\draw (6, 3.5) .. controls (6, 3) and (5, 2.5) .. (5, 1.5);
\draw (5, 1.5) .. controls (6, 0.75) and (6, 2.25) .. (7, 1.5);
\draw (8, 3.5) .. controls (8, 3) and (7, 2.5) .. (7, 1.5);
\draw (8, 3.5) .. controls (7, 4.25) and (7, 2.75) .. (6, 3.5);
\draw (-0.25, 4.55) node{{\tiny \(\proj^2\)}};
\draw (7, 4.75) node{{\tiny \(\proj^{23}\)}};
\draw[->] (0.8, 4.15) -- (5.755, 3);
\draw (1.6, 4.6) node{{\tiny \([x: y: z] \mapsto\)}};
\draw (3.0, 4.2) node{{\tiny \([x^7 + \cdots : \cdots : \cdots + 2z^7]\)}};
\draw (2.4, 5.0) node{{\tiny nonlinear map \(\proj^2 \to \proj^{23}\):}};
\draw (7.25, 3.5) .. controls (7.5, 2.5) and (7.2, 3) .. (6.5, 3);
\draw (6.5, 3) .. controls (5.5, 3) and (5.5, 2.5) .. (6.5, 2.5);
\draw (6.5, 2.5) .. controls (7, 2.5) and (7, 1.5) .. (6.5, 2);
\draw (6.5, 2) .. controls (5.75, 2.75) and (6.25, 1.75) .. (5.25, 1.75);
\draw[->, densely dotted] (-1, 2) -- (-0.5, 3.68);
\draw (-1, 1.8) node{{\tiny degree \(14\)}};
\draw (-1, 1.5) node{{\tiny curve in \(\proj^2\)}};
\draw (-1, 1.1) node{{\tiny \(x^{14} + \cdots = 0\)}};
\draw[->, densely dotted] (-0.15, 0.95) .. controls (0.85, -0.05) and (3.28, -0.8) .. (6.28, 2.2);
\draw (2.25, -0.15) node{{\tiny turns into degree \(7 \times 14 = 98\) curve in \(\proj^{23}\)}};
\end{tikzpicture}
\caption{A family of curves in $\proj^{23}$ obtained as curves in $\proj^2$, which in turn is mapped into $\proj^{23}$}
\label{f:crazycurve}
\end{figure}

Even though this first family is maximal, there is still a sense in which the curves in this first family are all ``non-general'', as we now explain.  Although degree $14$ plane curves have genus \(78\), almost all of the genus $78$ curves are not degree \(14\) plane curves.  (In fact, the space of genus $78$ curves forms a single irreducible space of dimension $3 \times 78 -3 = 231$;  the space of degree $14$ plane curves forms a single irreducible space of dimension only \(\binom {16} 2-1 - 8 = 111\).)      More formally:  the locus of genus $78$ curves that are plane curves form a proper closed locus of smaller dimension.   More informally:  if you were to choose a random genus $78$ curve, it would not appear in the first family.

\epoint{There is a second,  ``general'', family, and it does not behave pathologically}
 There is another family, of dimension $812$, that comes from the Brill--Noether theorem.  This second family is also maximal.
  The curves in the second family include ``almost all'' genus $78$ curves.  By this we mean, informally, that if you were to choose a random genus $78$ curve, it would appear in the second family.  More formally:  the locus of genus $78$ curves {\em not} appearing in the second family forms a proper closed locus of smaller dimension.
    It turns out that if we exclude components of $\cm(g,r,d)$ that consist of ``non-general curves" (curves from  some proper closed locus of smaller dimension in the moduli space of curves), then at most  one component remains, and we can know its dimension! 
    More precisely, let $\rho(g,r,d) = g - (r+1)(g-d + r)$.    
    
\tpoint{Brill--Noether Theorem}
{\em \label{t:BN}If $\rho \geq 0$, there is precisely one irreducible component of $\cm(g,r,d)$ consisting of general curves of genus \(g\), which we denote (for this article only) $\cm(g,r,d)^{BN}$.   It has dimension $\dim \cm(g,r,d)^{BN} = (r+1)d-(r-3)(g-1)$.
We call
such curves in $\proj^r$ \defi{Brill--Noether curves}.  If $\rho<0$, then there is no such component.  } \exercisedone

    This is the celebrated {\em Brill--Noether Theorem}. (We have given a very particular motivation for it, but there are many others.)   One  sign of an important theorem is that there are many different proofs giving many different insights, and this is one such example.    The original proof in the 1980's spans the four papers  \cite{E1, E2, E3, E4}; other alternative proofs from that era include \cite{E10, E11, E12}.  These proofs are valid in characteristic zero, but in fact the result is true in arbitrary characteristic, and that will be relevant for the upcoming discussion.  A proof valid in positive characteristic is included in the union of \cite{E5, E6, E7, E8, E9}.    Also, in the case when $\rho=0$, the fiber of \(\cm(g, r, d)\) over a general curve of genus \(g\) is a finite set of points.   We remark that this {\em number} is a beautiful combinatorial constant (the ``enumerative part'' of the Brill--Noether theorem); the first proofs of this were \cite{E13, E14}.  And as an aside, the failure of the theorem for curves with a fixed degree cover of $\proj^1$ (``gonality'') has led to some rich recent results, see for example \cite{E15, E16, JR, E17, E9, P}.
    
A little geometric intuition will suggest that if you have a family of curves in  $\proj^r$ of dimension $N$, then the locus of such curves passing through a fixed (generally chosen) point $p \in \proj^r$ should be empty or of dimension $N-(r-1)$.    The combination of that insight with the Brill--Noether Theorem~\ref{t:BN} 
(generalizing the discussion of \S \ref{s:countingpointconditions})
yields the following hope.

    \tpoint{Naive Conjecture}\label{c:naive}  {\em Brill--Noether curves of degree \(d\) and genus \(g\) in \(\proj^r\) interpolate
  \begin{equation}\label{eq:bn_conj}
        \left\lfloor\frac{\dim \cm(g,r,d)^{BN}   }{\dim \proj^r-1} \right\rfloor = \left\lfloor\frac{(r+1)d - (r-3)(g-1)}{r-1} \right\rfloor
  \end{equation}
 general points.}

 The expected number of general points in Conjecture~\ref{c:naive} is a priori an {\em upper bound} on the truth.  The reason is that, as we continue to impose the condition of interpolating each additional general point, the dimension drops by \(r-1\) each time, unless it prematurely becomes empty.

\epoint{Side Remark}
We remark that if you 
work with other types of curves besides Brill--Noether curves, then the expectation in Conjecture~\ref{c:naive} may not even provide an upper bound.  For example, 
Conjecture~\ref{c:naive} predicts that Brill--Noether curves of degree \(98\) and genus \(78\) in \(\proj^{23}\) interpolate \(36\) general points.  However, one can show that degree 98 genus 78 curves in \(\proj^{23}\) in the first family described in \S\ref{s:first_family} interpolate strictly more points.

  \section{Interpolation for Brill--Noether curves: the Main Theorem}

  Equipped with this Naive Conjecture~\ref{c:naive}, we can begin to explore the interpolation problem for Brill--Noether curves.  We begin with a series of examples.

  \bpoint{Example: Curves of degree \(4\) and genus \(1\) in \(\proj^3\)}\label{ex:deg4gen1}\label{ex:d4g1}
Conjecture~\ref{c:naive} predicts that such curves can interpolate \(8\) general points.

  This is one case in which the classical approach for plane curves can be successfully employed for Brill--Noether curves in higher dimensions.  What makes this case simpler is that every (Brill--Noether) curve of degree \(4\) and genus \(1\) in \(\proj^3\) can be described as the simultaneous vanishing of two quadratic equations.  Geometrically, we say that these curves are the complete intersection of two quadric surfaces.  (For readers familiar with algebraic curves or Riemann surfaces, this is a consequence of the Riemann--Roch theorem.)

  As we saw in \S\ref{s:quadric_surface}, the parameter space of quadric surfaces is \(\proj^9\), and those quadric surfaces passing through a given point form a codimension \(1\) linear space in \(\proj^9\).  Intersecting \(8\) of these codimension \(1\) linear spaces leaves a \(1\)-dimensional linear space, which is the span of two independent quadratic equations \(Q_1, Q_2\).  By construction, these surfaces each pass through the \(8\) general points, and so the degree \(4\) and genus \(1\) curve \(Q_1 \cap Q_2\) does as well.  This proves the Conjecture~\ref{c:naive} in this case (and in fact, shows that there is a {\em unique} Brill--Noether curve of degree \(4\) and genus \(1\) through \(8\) general points).

\bpoint{Example: Curves of degree \(5\) and genus \(2\) in \(\proj^3\)}\label{ex:d5g2} Conjecture~\ref{c:naive} predicts that such curves can interpolate \(10\) general points.

Similarly to Example~\ref{ex:deg4gen1}, every such curve lies on a (unique) quadric surface (see the next paragraph for more details). But as opposed to the case of degree \(4\) and genus \(1\) in \S\ref{ex:d4g1}, the quadric surface containing \(C\) gets in the way!  If Conjecture~\ref{c:naive} is true for \((g,r,d) = (2, 3, 5)\), then given 10 general points \(p_1, \dots, p_{10}\) in \(\proj^3\), there exists some (Brill--Noether) curve \(C\) of degree \(5\) and genus \(2\) passing through these points.  Since \(C\) is always contained in a quadric surface \(Q\), and such a surface containing \(C\) also contains \(p_1, \dots, p_{10}\), quadric surfaces would interpolate \(10\) general points.  But as we saw in \S\ref{s:quadric_surface}, linear algebraic methods prove that quadric surfaces interpolate only \(9\) general points!  What we have on our hands is a counterexample to Conjecture~\ref{c:naive}!

We can explain the existence of the quadric surface in two ways.
First, it is a consequence of the Riemann--Roch theorem.  We have already observed that there is a \(10\)-dimensional vector space of quadratic equations on \(\proj^3\) (giving a \(9\)-dimensional projective space \(\proj^9\) of quadric surfaces).  In order for a quadric surface to contain a curve \(C\), the restriction of (any) equation must vanish {\em identically} along the curve.  Restriction to \(C\) is a linear map of vector spaces, and by the Riemann--Roch theorem, the codomain (quadratic equations on \(C\)) is a vector space of dimension \(9\).  Consequently, this map {\em must} have a nontrivial kernel, which gives a quadric surface containing \(C\).
Since it will be relevant in other (counter)examples, we can also give a synthetic construction of the quadric surface \(Q \supset C\), relying on a bit more geometry of curves.  A curve of genus \(2\) is always \defi{hyperelliptic}, i.e., it admits a \(2\)-to-\(1\) map \(\pi \colon C \to \proj^1\).  
The quadric surface \(Q\) is obtained as the union of all lines spanned by pairs of points \(p\) and \(q\) for which \(\pi(p) = \pi(q)\).  A surface obtained in this way we call a \defi{hyperelliptic scroll}. 

It turns out that degree \(5\) and genus \(2\) in \(\proj^3\) is only the first in a (finite) series of counterexamples that all have a similar flavor.  These counterexamples all arise because the given type of curve is contained in a certain type of surface, and that surface fails to interpolate the required number of general points.  In the case \((g, r, d) = (2, 3, 5)\), the surface was a quadric surface (that also happened to be a hyperelliptic scroll).

\bpoint{Further counterexamples: obstructions from surfaces}\label{s:counter}

\epoint{(Counter-)Example: Curves of degree \(7\) and genus \(2\) in \(\proj^5\)}\label{ex:d7g2} Conjecture~\ref{c:naive} predicts that such curves can interpolate \(10\) general points. Since \(C\) has genus \(2\), it is hyperelliptic and so again lies on a hyperelliptic scroll \(S \subset \proj^5\) (which turns out to be a scroll of degree \(4\)).  Such scrolls interpolate 9 general points, and so again force Conjecture~\ref{c:naive} to be false in this case.

\epoint{(Counter-)Example: Canonical curves in \(\proj^3\) (degree \(6\) and genus \(4\))}  
When \((g,r,d) = (g, g-1, 2g-2)\), the realization of the curve in projective space is canonically attached to the curve (and is called a \defi{canonical curve of genus \(g\)}).  Since it is canonical, the geometry of this projective embedding reflects deep intrinsic geometry of the curve.  For this reason, canonical curves have been intensely studied for many years and are the subject of many important results and conjectures.

Conjecture~\ref{c:naive} predicts that canonical curves in \(\proj^3\) interpolate \(12\) general points.  The Riemann--Roch theorem can again be used to prove that any such curve is specified by the simultaneous vanishing of a quadratic equation and a cubic equation, i.e., it is a complete intersection of a quadric surface and cubic surface.  Consequently, any such curve is again contained in a quadric surface, which only interpolates 9 general points, instead of 12.

\epoint{(Counter-)Example: Canonical curves in \(\proj^5\) (degree \(10\) and genus \(6\))} \label{ex:d10g6} 
Conjecture~\ref{c:naive} predicts that canonical curves in \(\proj^5\) interpolate \(12\) general points. An explicit description of (general) canonical curves of genus \(6\) is the intersection of a quadric hypersurface and a del Pezzo surface of degree \(5\).
In this case, the del Pezzo surface can only interpolate 11 general points, and so again obstructs Conjecture~\ref{c:naive} in this case.

\bpoint{Where do all of these counterexamples leave us?}
It is easy to be discouraged by the presence of the four counterexamples enumerated in \S\ref{ex:d5g2}--\ref{ex:d10g6}: these are some of the first examples of Brill--Noether curves you might investigate (low genus curves and canonical curves) and already they do not behave as the naive dimension estimates of Conjecture~\ref{c:naive} predict!  It is worth noting, however, that these counterexamples naturally lie in two infinite families (curves of genus \(2\) and degree \(r+2\) in \(\proj^r\), and canonical curves), and the reason that these examples are counterexamples do not persist to the rest of the family.  This is one reason for optimism that these counterexamples are coincidences of small numbers, and are not indicative of the general behavior.

In fact, the main theorem of \cite{lv} says precisely this: these four examples in \S\ref{ex:d5g2}--\ref{ex:d10g6} are the \emph{only} counterexamples to Conjecture~\ref{c:naive} for Brill--Noether curves.

  \tpoint{Main Theorem (Interpolation for Brill--Noether curves) \cite{lv}}\label{t:main}  {\em A Brill--Noether curve of degree \(d\) and genus \(g\) in \(\proj^r\) interpolates the expected number
  \[
  \left\lfloor\frac{(r+1)d - (r-3)(g-1)}{r-1} \right\rfloor
  \]
  of general points if and only if
  \[
  (g, r, d) \not\in \{ (2, 3, 5), (4, 3, 6), (2, 5, 7), (6, 5, 10) \}.
  \]}

  The ``only if'' portion of Theorem~\ref{t:main} is straightforward, as we have sketched in \S\ref{ex:d5g2}--\ref{ex:d10g6}.  The difficult and novel aspect is the proof that there are no other counterexamples.

  \bpoint{High-level overview of the proof strategy of Theorem~\ref{t:main}}
  Theorem~\ref{t:main} is really a triply-infinite collection of statements, one for each tuple \((g, r, d)\) such that \(\rho(g, r, d) \geq 0\).  At a very high level, the proof proceeds by induction, reducing one case of Theorem~\ref{t:main} to another with smaller values of \((g, r, d)\).  Since the values \((g, r, d)\) have geometric significance as the invariants of a Brill--Noether curve, you might guess that this reduction is achieved by {\em breaking the curve} into pieces with smaller invariants.  This procedure is delicate for a number of reasons, the most obvious being the presence of counterexamples!  Any inductive argument to reduce \((g, r, d)\) {\em must} break down before hitting one of \(\{ (2, 3, 5), (4, 3, 6), (2, 5, 7), (6, 5, 10) \}\).  In actuality, there are multiple inductive steps that combine to reduce all cases of Theorem~\ref{t:main} to about \(30\) base cases that form a ``shield'' around the counterexamples at the bottom.  

  To illustrate the potential of breaking a curve into simpler pieces, we now give  another (third!) proof that conics interpolate 5 general points by this method.
  A more detailed and technical explanation of the inductive steps appears in \S\ref{s:arch}.
  
\bpoint{(A third) Proof that conics interpolate 5 general points}\label{s:conics2}
It is easy to make conics that pass through \(5\) (potentially special) points: just choose any \(5\) points on any conic.
If conics interpolate fewer than \(5\) {\em general} points, then the collections of \(5\) points that {\em do} lie on a conic must be contained in the vanishing locus of some (nonzero) equation in \((\proj^2)^5\).  To rule this out, it suffices to show that given some collection of \(5\) points on some conic, we can ``deform'' these points in {\em any} small direction while preserving the property that they lie on a conic.  

Since we said ``some'' collection of \(5\) points on ``some'' conic, we are free to choose these to be quite degenerate!  Namely, we will take the conic to be the union of the two coordinate axes, which is specified by the quadratic equation \(xy = 0\).  Any collection of \(5\) points \(p_1, \dots, p_5\) on this degenerate conic must contain \(3\) collinear points, and hence look quite special.  Suppose that we begin with \(3\) points on the \(x\)-axis and \(2\) points on the \(y\)-axis, as pictured on the left in Figure~\ref{fig:broken_conic}.  Consider an arbitrary deformation of these points; that is, we ``bump'' each of these points in some arbitrarily small direction.  We can formalize this by adding some small quantity \((\alpha_i \epsilon, \beta_i\epsilon)\) to each point \(p_i\), where \(\epsilon>0\) is an infinitesimal.  These deformed points are illustrated with unfilled circles  in Figure~\ref{fig:broken_conic}.  Our goal is to similarly ``deform'' the degenerate conic \(xy=0\) to continue to pass through the \(5\) points.  An infinitesimal deformation of \(xy=0\) will be specified by the vanishing of an equation of the form \(xy + \epsilon Q(x,y)\) for some quadratic equation \(Q(x, y)\).

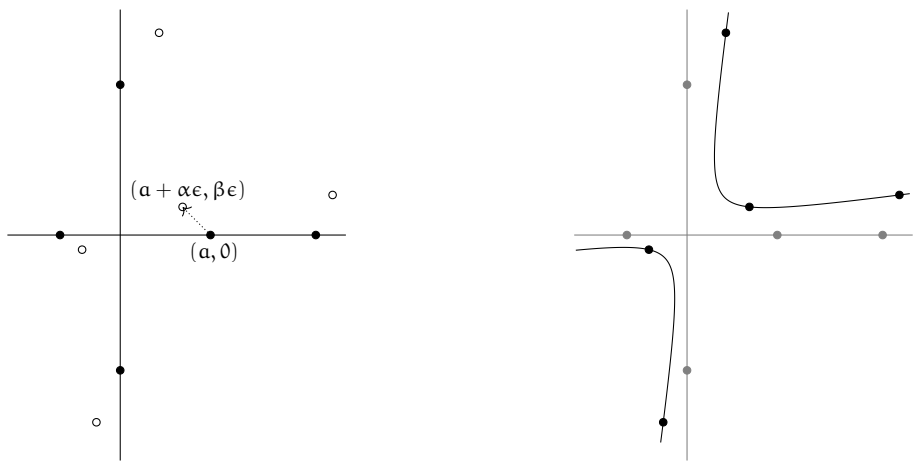
\begin{figure}[h!]
    \centering
\begin{minipage}{.45\textwidth}
\begin{center}
\begin{tikzpicture}
\draw (0,-3) -- (0, 3);
\draw (-1.5,0) -- (3,0);
\filldraw (-.8,0) circle[radius=0.05];
\filldraw (1.2,0) circle[radius=0.05];
\filldraw (2.6,0) circle[radius=0.05];
\filldraw (0,2) circle[radius=0.05];
\filldraw (0, -1.8) circle[radius=0.05];
\coordinate (A) at (.6,.6);
\coordinate (B) at (-.4,-.4);
\pgfmathsetmacro{\acRatio}{0.8}

\coordinate (BA) at ($ (B)-(A) $);
\newdimen\myBAx
\pgfextractx{\myBAx}{\pgfpointanchor{BA}{center}}
\newdimen\myBAy
\pgfextracty{\myBAy}{\pgfpointanchor{BA}{center}}
\pgfmathsetlengthmacro{\c}{veclen(\myBAx,\myBAy)/2}
\pgfmathsetlengthmacro{\b}{sqrt(1-\acRatio^2)*\c}
\pgfmathsetlengthmacro{\a}{\acRatio*\c}
\pgfmathanglebetweenlines{\pgfpoint{0}{0}}{\pgfpoint{1}{0}}
{\pgfpointanchor{A}{center}}{\pgfpointanchor{B}{center}}
\let\rotAngle\pgfmathresult
\coordinate (O) at ($ (A)!.5!(B) $);
\tikzset{hyperbola1/.style={rotate=\rotAngle,shift=(O),
    domain=-2.1:2.1,variable=\t,samples=50,smooth}}
\tikzset{hyperbola2/.style={rotate=\rotAngle,shift=(O),
    domain=-1.5:2.1,variable=\t,samples=50,smooth}}
\draw[hyperbola2, color=white] plot ({ \a*cosh(\t)},{\b*sinh(\t)});
\draw[hyperbola1, color=white] plot ({-\a*cosh(\t)},{\b*sinh(\t)});

\draw[hyperbola2, color=black] ({ -\a*cosh(-2)},{\b*sinh(-2)}) circle[radius=0.05];
\draw[hyperbola2, color=black] ({ -\a*cosh(.7)},{\b*sinh(.7)}) circle[radius=0.05];
\draw[hyperbola2, color=black] ({ -\a*cosh(2.05)},{\b*sinh(2.05)}) circle[radius=0.05];

\draw[hyperbola1, color=black] ({ \a*cosh(-.5)},{\b*sinh(-.5)}) circle[radius=0.05];
\draw[hyperbola1, color=black] ({ \a*cosh(2)},{\b*sinh(2)}) circle[radius=0.05];

\coordinate (S) at (.83, .37); 

\draw (0.9, .3) node[above] {{\tiny \color{black}\((a + \alpha \epsilon, \beta \epsilon)\)}};
\draw (1.25, 0.05) node[below] {{\tiny \((a,0)\)}};
\draw[color=black, ->, densely dotted] (1.2,0) -- (S);

\end{tikzpicture}
\end{center}
\end{minipage}
\begin{minipage}{.45\textwidth}
\begin{center}
\begin{tikzpicture}
\draw[color=gray] (0,-3) -- (0, 3);
\draw[color=gray] (-1.5,0) -- (3,0);

\filldraw[color=gray] (-.8,0) circle[radius=0.05];
\filldraw[color=gray] (1.2,0) circle[radius=0.05];
\filldraw[color=gray] (2.6,0) circle[radius=0.05];
\filldraw[color=gray] (0,2) circle[radius=0.05];
\filldraw[color=gray] (0, -1.8) circle[radius=0.05];
\coordinate (A) at (.6,.6);
\coordinate (B) at (-.4,-.4);
\pgfmathsetmacro{\acRatio}{0.8}

\coordinate (BA) at ($ (B)-(A) $);
\newdimen\myBAx
\pgfextractx{\myBAx}{\pgfpointanchor{BA}{center}}
\newdimen\myBAy
\pgfextracty{\myBAy}{\pgfpointanchor{BA}{center}}
\pgfmathsetlengthmacro{\c}{veclen(\myBAx,\myBAy)/2}
\pgfmathsetlengthmacro{\b}{sqrt(1-\acRatio^2)*\c}
\pgfmathsetlengthmacro{\a}{\acRatio*\c}
\pgfmathanglebetweenlines{\pgfpoint{0}{0}}{\pgfpoint{1}{0}}
{\pgfpointanchor{A}{center}}{\pgfpointanchor{B}{center}}
\let\rotAngle\pgfmathresult
\coordinate (O) at ($ (A)!.5!(B) $);
\tikzset{hyperbola1/.style={rotate=\rotAngle,shift=(O),
    domain=-2.1:2.1,variable=\t,samples=50,smooth}}
\tikzset{hyperbola2/.style={rotate=\rotAngle,shift=(O),
    domain=-1.5:2.1,variable=\t,samples=50,smooth}}
\draw[hyperbola2, color=black] plot ({ \a*cosh(\t)},{\b*sinh(\t)});
\draw[hyperbola1, color=black] plot ({-\a*cosh(\t)},{\b*sinh(\t)});

\filldraw[hyperbola2, color=black] ({ -\a*cosh(-2)},{\b*sinh(-2)}) circle[radius=0.05];
\filldraw[hyperbola2, color=black] ({ -\a*cosh(.7)},{\b*sinh(.7)}) circle[radius=0.05];
\filldraw[hyperbola2, color=black] ({ -\a*cosh(2.05)},{\b*sinh(2.05)}) circle[radius=0.05];

\filldraw[hyperbola1, color=black] ({ \a*cosh(-.5)},{\b*sinh(-.5)}) circle[radius=0.05];
\filldraw[hyperbola1, color=black] ({ \a*cosh(2)},{\b*sinh(2)}) circle[radius=0.05];

\coordinate (S) at (.83, .37); 

\end{tikzpicture}
\end{center}
\end{minipage}
    \caption{Deformation strategy to prove that conics interpolate \(5\) general points.}
    \label{fig:broken_conic}
\end{figure}

Let us focus on one particular original point \(p = (a,0)\) on the \(x\)-axis, and suppose that we deformed it to the point \((a + \alpha\epsilon, \beta\epsilon)\).  In order for our deformed conic to pass through this point, we need
\((a + \alpha \epsilon)(\beta \epsilon) + \epsilon Q(a + \alpha \epsilon, \beta\epsilon) = 0\).
By a limit procedure, it suffices to ignore all of the \(\epsilon^2\) terms in this equation, and only consider the terms that are linear in our infinitesimal \(\epsilon\).  This substantially simplifies the equation to \(Q(a, 0) = -a\beta\).  The exact equation is not as relevant as the fact this the condition of the {\em deformed} conic passing through the {\em deformed} point is equivalent to the quadratic equation \(Q(x, 0)\) in \(x\) taking a specified value at the {\em original} point.  We have three such points, and Lagrange interpolation guarantees that we can find a quadratic equation \(Q_1(x,0)\) taking the required values!  This takes care of the interpolation condition at all of the points that are deformations of points on the \(x\)-axis.  We can then find a quadratic equation \(Q_2(0, y)\) in \(y\) taking the required values at the two original points on the \(y\)-axis.  (The required values here depend on both the deformations of the two points on the \(y\)-axis and the constant term \(Q_1(0,0)\) previously found.)  It might look like we have another degree of freedom from Lagrange interpolation, but we must also require that \(Q_2(0,0) = 0\) so that we don't disturb the constant term of \(Q_1(x,0)\) when we set \(Q(x,y) = Q_1(x,0) + Q_2(0, y)\).

This exact argument doesn't generalize to Brill--Noether curves in higher-dimensional projective space, since we relied on the equation of the curve to understand how to deform it.  Nevertheless, degeneration to broken curves combined with deformation theory is a powerful tool to understand the interpolation problem.  

The key insight that made this proof work was the fact that deformations of the union of the two axes were flexible enough to be able to specify where \(5\) points go in the deformation.
Outside of the plane, deformations are no longer controlled by an equation, but rather by the \defi{normal bundle} of the curve.  This is the topic of the next section.

\bpoint{Normal bundles and deformation theory}
The \defi{normal space} of a curve \(C \subset \proj^r\) at a point \(p \in C\) is the quotient \(T_p \proj^r/T_pC\) of tangent directions in \(\proj^r\) at \(p\) by the tangent direction along \(C\) at \(p\).  These \((r-1)\)-dimensional vector spaces fit together as the point \(p\) varies into the \defi{normal bundle} \(N_C\) of the curve.  (Even though this definition only makes sense for smooth embedded curves \(C \subset \proj^r\), it can be extended to curves with mild singularities --- in particular nodal singularities where two branches meet transversely.  At the singular points, the fiber of the normal bundle is {\em not} simply \(T_p \proj^r/T_pC\), which is one of the subtleties of this strategy.)

\begin{figure}[h!]
\begin{tikzpicture}[scale=1]
\draw (0, 0) .. controls (1, 3) and (2, 3) .. (3, 1);
\draw (3, 1) .. controls (4, -1) and (4, 2) .. (5, 3);
\draw[densely dotted] (-0.3, 0.3) .. controls (1, 3) and (2, 2.5) .. (3, 1);
\draw[densely dotted] (3, 1) .. controls (4, -0.5) and (4, 2) .. (5, 2.5);
\draw[->] (0.15, 0.43) -- (-0.17, 0.56);
\draw[->] (0.3, 0.82) -- (0.03, 0.93);
\draw[->] (0.45, 1.15) -- (0.24, 1.27);
\draw[->] (0.6, 1.45) -- (0.44, 1.55);
\draw[->] (0.8, 1.78) -- (0.7, 1.84);
\draw[->] (1.5, 2.373) -- (1.53, 2.224);
\draw[->] (1.9, 2.34) -- (1.83, 2.16);
\draw[->] (2.225, 2.13) -- (2.11, 2.005);
\draw[->] (2.5, 1.83) -- (2.41, 1.75);
\draw[->] (3.6, 0.33) -- (3.58, 0.53);
\draw[->] (3.84, 0.54) -- (3.74, 0.61);
\draw[->] (3.31, 0.5) -- (3.39, 0.575);
\draw[->] (4.495, 2.2) -- (4.6, 2.15);
\draw[->] (4.64, 2.5) -- (4.845, 2.4);
\end{tikzpicture}
\caption{``deformation of a curve" $\sim$ ``section of its normal bundle"}
\label{fig:section_normal}
\end{figure}
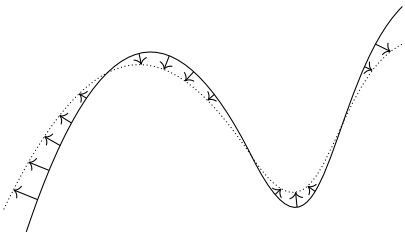

An infinitesimal deformation of \(C \subset \proj^r\) is specified by giving a normal vector at each point of the curve, which labels ``where that point goes in the deformation'' as in Figure~\ref{fig:section_normal}.  Such a choice of normal vector at every point of \(C\) is called a \defi{section} of the normal bundle of \(C\).  We can formalize the property that ``deformations are flexible enough to specify where the appropriate number of general points deform to'' in terms of sections of the normal bundle.  Namely, we want to be able to find a section of the normal bundle taking specified values at each of the \(n\) points.  In other words, we want the evaluation map
\[\{\text{sections of \(N_C\)}\} \to \bigoplus_i N_C|_{p_i},\]
taking a section to its value at each \(p_i\), to be surjective (when it should be).  This is equivalent to requiring that the kernel --- those sections that vanish at \(p_1, \dots, p_n\) --- has the correct dimension.  This leads to the following definition.

\bpoint{Definition: Interpolation for the normal bundle} 
We say that the normal bundle \(N_C\) of a Brill--Noether curve \(C \subset \proj^r\) satisfies \defi{interpolation} if, for all \(n \geq 0\) and for general points \(p_1, \dots, p_n\) on \(C\), the dimension of the space of sections of \(N_C\) that {\em vanish} at \(p_1, \dots, p_n\) is exactly
\[\max\big(0, (r+1)d - (r-3)(g-1) - (r-1)n\big).\]

An argument in deformation theory that is a souped-up version of \S\ref{s:conics2} shows that if there exists some (possibly degenerate) Brill--Noether curve \(C \subset \proj^r\) of degree \(d\) and genus \(g\) for which \(N_C\) satisfies interpolation, then Conjecture~\ref{c:naive} holds for \((g, r, d)\).  

The converse is {\em not} true!  This property of the normal bundle can be strictly stronger than Conjecture~\ref{c:naive}; in fact, it implies a version of Conjecture~\ref{c:naive} where the curve can be chosen to meet linear spaces of arbitrary dimension (the dimension \(0\) case is points).  This leads to an additional counterexample from the family of \S\ref{ex:d5g2} and \S\ref{ex:d7g2}.  Curves of degree \(6\) and genus \(2\) are expected to be able to pass through \(9\) general points {\em and} meet a general line; they can pass through \(9\) general points, but the hyperelliptic scroll containing them prevents them from simultaneously passing through \(9\) general points and meeting a general line.  Hence the normal bundle of such a curve never satisfies interpolation.

\bpoint{Exotic counterexamples in characteristic \(2\)}
Something else can go wrong over fields of characteristic \(2\).
To see this, let's return to the definition of the normal bundle as the quotient \(T_{\proj^r}|_C/T_C\) of the tangent bundle to projective space along \(C\) by the tangent bundle of \(C\).  How does \(T_C\) sit inside \(T_{\proj^r}|_C\)?  Explicitly, if \(C = \proj^1\), so that a map \(C \to \proj^r\) corresponds to an \((r+1)\)-tuple of degree \(d\) polynomials, then the partial derivatives of these polynomials give the inclusion \(T_C \to T_{\proj^r}|_C\).  This is the familiar definition of a derivative: if I move a tiny bit in some direction in \(C\), then the derivative(s) of \(f\) measures how much the image moves in \(f(C)\).

The derivative has a funny behavior in characteristic \(2\): the derivative of any monomial \(x^a\) for even \(a\) is zero.  This means that the matrix of partial derivatives of a tuple of polynomials \((f_0, \dots, f_r)\) does not look general: it only has even powers!  This manifests in a parity constraint on ``the degrees of polynomials defining the normal bundle'' of a curve of genus \(0\) in projective space.  This parity constraint means that \(N_C\) can only satisfy interpolation if
the quantity appearing inside the floor symbols in \eqref{eq:bn_conj} is an integer with the opposite parity as \(d\), or explicitly, 
\(d \equiv 1 \pmod{r-1}\).   In other words:  interpolation  (for degree $d$ genus $0$ curves in $\proj^r$ in characteristic $2$)  is {\em guaranteed} to fail if $d \not\equiv 1 \pmod {r-1}$.

For experts, we cast this in more technical language.  In characteristic \(2\), the twist of the normal bundle \(N_C(-1)\) of a curve \(C\) is isomorphic to the pullback under the Frobenius morphism \(F \colon C \to C'\) of a bundle \(N'\) on the Frobenius twist \(C'\) of \(C\).
This is not usually a problem, but if \(C\) has genus \(0\), then since \(F\) has degree \(2\), the normal bundle \(N_C\) must necessarily be isomorphic to a direct sum of line bundles whose degree has the {\em same} parity as the degree \(d\) of \(C\).
This translates into the stated congruence constraint.

Here is why you should care even if you are not an expert, or don't care about positive characteristic:  If your arguments are fully algebraic, and they don't somehow bring in the characteristic, then they have to fail!  Our ``guardrails" in our induction have  to somehow protect us from such craziness.  From another point of view, this characteristic 2 pathology  tells us what approaches cannot work in characteristic $0$.

\tpoint{Theorem (Interpolation for normal bundles of Brill--Noether curves) \cite{lv}}\label{t:main_normal} {\em The normal bundle of a general Brill--Noether curve of degree \(d\) and genus \(g\) in \(\proj^r\) satisfies interpolation if and only if
\newline \indent $\bullet$  \((g, r, d) \not\in \{ (2, 3, 5), (4, 3, 6), (2, 4, 6), (2, 5, 7), (6, 5, 10) \}\), and
\newline \indent $\bullet$ 
   if the characteristic is \(2\) and \(g=0\), then \(d \equiv 1 \pmod{r-1}\).}

The additional characteristic \(2\) counterexamples in Theorem~\ref{t:main_normal} do not cause problems for Conjecture~\ref{c:naive}, even in characteristic \(2\), but an additional argument is needed in this case to deduce these cases of Theorem~\ref{t:main}.
Roughly speaking, we take the coordinates of the points through which we want to find a curve in characteristic \(2\) and lift them to characteristic \(0\).  By the interpolation result in characteristic \(0\), we can find a suitable curve interpolating the lifted points.  We finally reduce the defining equation of this curve modulo \(2\) to obtain an interpolating curve in characteristic \(2\).  What could go wrong?  The subtle point is whether the reduced curve is a Brill--Noether curve of the same degree and genus as the curve in characteristic \(0\).  Fortunately, 
the only additional exceptions in characteristic \(2\) occur when \(g=0\), and it is easy to understand the possible types of the reduction of a curve of genus \(0\).
  
\section{Architecture of the proof}\label{s:arch}

We now get into some details.

The proof of Theorem~\ref{t:main_normal} is by induction.
More precisely, there are two basic inductive strategies:
breaking the curve into simpler pieces in a fixed projective space,
and exploiting projection maps to reduce to projective spaces of smaller dimension.

\bpoint{Breaking the curve\label{ss:break}}
We specialize the curve \(C\) to a reducible curve \(X \cup Y\).

The first subtlety here is that the restriction \(N_{X \cup Y}|_X\) does
\emph{not} coincide with \(N_X\); something interesting happens at \(X \cap Y\).
For example, suppose \(X \cup Y\) is the union of the two coordinate axes in the plane.
Then, a smoothing of \(X \cup Y\) to a hyperbola, interpreted as a section of the normal bundle of the \(x\)-axis,
has a pole at the origin, as illustrated in Figure~\ref{fig:hyperbola}.

\begin{figure}[h!]
\begin{tikzpicture}
\draw (-3, 0) -- (3, 0);
\draw (0, -3) -- (0, 3);
\draw[densely dotted] (0.12, 3) .. controls (0.28, 0.28) .. (3, 0.12);
\draw[densely dotted] (-0.12, -3) .. controls (-0.28, -0.28) .. (-3, -0.12);
\draw[->] (2.75, 0) -- (2.75, 0.1309090909090909);
\draw[->] (2.5, 0) -- (2.5, 0.144);
\draw[->] (2.25, 0) -- (2.25, 0.16);
\draw[->] (2, 0) -- (2, 0.18);
\draw[->] (1.75, 0) -- (1.75, 0.20571428571428574);
\draw[->] (1.5, 0) -- (1.5, 0.24);
\draw[->] (1.25, 0) -- (1.25, 0.288);
\draw[->] (1, 0) -- (1, 0.36);
\draw[->] (0.75, 0) -- (0.75, 0.48);
\draw[->] (0.5, 0) -- (0.5, 0.72);
\draw[->] (0.25, 0) -- (0.25, 1.44);
\draw[->] (-2.75, 0) -- (-2.75, -0.1309090909090909);
\draw[->] (-2.5, 0) -- (-2.5, -0.144);
\draw[->] (-2.25, 0) -- (-2.25, -0.16);
\draw[->] (-2, 0) -- (-2, -0.18);
\draw[->] (-1.75, 0) -- (-1.75, -0.20571428571428574);
\draw[->] (-1.5, 0) -- (-1.5, -0.24);
\draw[->] (-1.25, 0) -- (-1.25, -0.288);
\draw[->] (-1, 0) -- (-1, -0.36);
\draw[->] (-0.75, 0) -- (-0.75, -0.48);
\draw[->] (-0.5, 0) -- (-0.5, -0.72);
\draw[->] (-0.25, 0) -- (-0.25, -1.44);
\end{tikzpicture}
\caption{Hyperbola as a deformation of the union of coordinate axes}
\label{fig:hyperbola}
\end{figure}
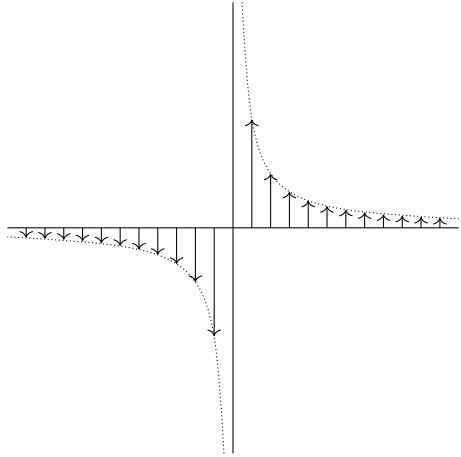

The story in higher dimensions is similar (see \cite{hh}): Sections of the restriction \(N_{X \cup Y}|_X\)
can be interpreted as meromorphic sections of \(N_X\), with simple poles allowed at \(X \cap Y\)
that point in the direction of \(Y\).
This procedure for obtaining \(N_{X \cup Y}|_X\) from \(N_X\) is known as a \defi{(positive) elementary modification}.

The second subtlety here is that interpolation for \(N_{X \cup Y}\) does not follow from
interpolation for \(N_{X \cup Y}|_X\) and \(N_{X \cup Y}|_Y\).
A vector bundle on a reducible nodal curve \(X \cup Y\) is more data than its restriction to each of the components.  It also records the data of the how the fibers are glued together over the nodes.  This gluing data is subtle and explains why interpolation for \(N_{X \cup Y}\) is not equivalent to interpolation for the two restrictions to \(X\) and to \(Y\).  Nevertheless, if one of the components \(Y\) is geometrically simple, we can get a geometric handle on the subtle gluing data.  Armed with this description, it is possible to encompass the gluing data into a modification of \(N_{X \cup Y}|_X\).  This is the machine powering our induction: we reduce from interpolation for \(N_{X \cup Y}\) to interpolation for a specific modification of \(N_{X \cup Y}|_X\). The rest of this subsection is fleshing out this idea in the necessary cases.

Unlike the first subtlety, there is no general method to circumvent this issue.
Nevertheless, for certain sufficiently simple curves \(Y\),
there are geometric constructions that allow us to relate interpolation for \(N_{X \cup Y}\) 
to interpolation for a \emph{further modification} of \(N_{X \cup Y}|_X\).
For example, suppose that \(Y = L\) is a line, meeting \(X\) at a single point \(u\).
Then \(N_{X \cup L}|_L\) is obtained from \(N_L \cong \mathcal{O}(1)^{\oplus (r - 1)}\)
by making a single modification at \(u\), so \(N_{X \cup L}|_L \cong \mathcal{O}(2) \oplus \mathcal{O}(1)^{\oplus(r - 2)}\).
In this case, interpolation for \(N_{X \cup L}\) turns out to follow from interpolation for
the modification of \(N_{X \cup L}|_X\) towards the fiber of the positive summand \(\mathcal{O}(2)|_u\).
Luckily, this admits a concrete geometric description:
The positive subbundle \(\mathcal{O}(2) \subset N_{X \cup L}|_L\)
is (the twist of) the normal bundle of \(L\) in the plane spanned by \(L\) and the tangent line \(T_u X\),
which is tangent along the entire line \(L\) to the cone over \(X\) with vertex \(v \in L \smallsetminus \{u\}\), as illustrated in Figure~\ref{fig:cone}.

\begin{figure}[h!]
\begin{tikzpicture}
\draw (2, 5) .. controls (0, 5) and (-1, 2.5) .. (0, 1.5);
\draw (0, 1.5) .. controls (0.5, 1) and (2, 0.75) .. (2, 1.5);
\draw (0.05, 1.55) .. controls (0.5, 2) and (2, 2.25) .. (2, 1.5);
\draw (-0.05, 1.45) .. controls (-1, 0.5) and (1, 0) .. (2, 0);
\filldraw (3, 4.125) circle[radius=0.05];
\draw[densely dotted] (0, 4.5) -- (4, 4);
\draw (0.5, 4.4375) -- (4, 4);
\draw[densely dotted] (3, 4.125) -- (-0.3, 3.75);
\draw[densely dotted] (3, 4.125) -- (-0.25, 4);
\draw[densely dotted] (3, 4.125) -- (-0.15, 4.25);
\draw[densely dotted] (3, 4.125) -- (0.2, 4.75);
\draw[densely dotted] (3, 4.125) -- (0.5, 5);
\draw[densely dotted] (3, 4.125) -- (0.9, 5.15);
\draw (2.1, 0) node{{\tiny \(X\)}};
\draw (3.5, 3.93) node{{\tiny \(L\)}};
\draw (3, 4) node{{\tiny \(v\)}};
\draw (0.45, 4.55) node{{\tiny \(u\)}};
\filldraw (0.5, 4.4375) circle[radius=0.05];
\draw[densely dotted] (1.5, 5.32) -- (-0.5, 3.56) -- (3, 3.1225) -- (5, 4.8825) -- (1.5, 5.32);
\end{tikzpicture}
\caption{Plane \(\langle L, T_u X \rangle\) is tangent along \(L\) to cone over \(X\) with vertex \(v\)}
\label{fig:cone}
\end{figure}
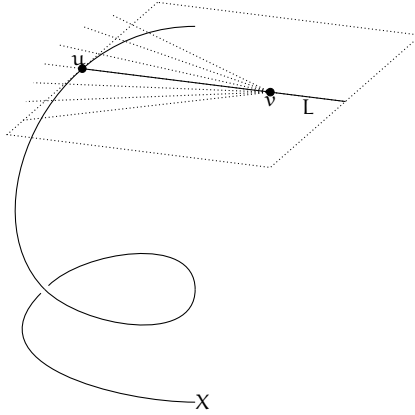

In particular, writing \(N_{X \to v}\) for the subbundle corresponding to ``normal directions that point towards \(v\)'',
i.e., for the normal bundle of \(X\) in this cone,
interpolation for \(N_{X \cup L}\) follows from a ``double modification'' of \(N_X\)
at \(u\) towards \(N_{X \to v}\).  The first modification takes us from \(N_X\) to \(N_{X \cup L}|_X\), and the second modification describes those sections that agree with the positive part of \(N_{X \cup L}|_L\).
We write this symbolically as \(N_X[2u \to v]\).

A similar formula can be derived for the case when \(Y = L\) is a \(2\)-secant line,
again using a similar geometric construction involving cones.
However, neither of these degenerations can increase the value of \(g + r - d\),
which can be arbitrarily large for a Brill--Noether curve.
Fortunately, when \(Y\) is a rational curve of degree \(r - 1\) meeting \(X\) in \(r + 1\) points,
the degree of \(N_{X \cup Y}|_Y\) is a multiple of \(r - 1 = \operatorname{rk} N_{X \cup Y}|_Y\).
In this case, one shows that \(N_{X \cup Y}|_Y\) is usually perfectly balanced,
and so we can reduce interpolation for \(N_{X \cup Y}\) to interpolation for \(N_{X \cup Y}|_X\).
These three degenerations thus form the backbone of the first inductive strategy.

The main defect of this type of inductive strategy is that each time we apply it,
we increase the number of modifications that are made to the normal bundle.
It is not good to have an inductive argument where something gets worse!

\bpoint{Projection maps\label{ss:proj}}
For any point \(p\), we have a projection map \(\mathbb{P}^r \dashrightarrow \mathbb{P}^{r - 1}\).
Then, writing \(\bar{C}\) for the projection of \(C\) from a point \(p \in C\), we have an exact sequence
\[0 \to N_{C \to p} \to N_C \to N_{\bar{C}}(p) \to 0.\]
Similarly, if we make modifications at \(x_1 + \cdots + x_n\) pointing towards \(p\), we obtain the sequence
\begin{equation} \label{proj}
0 \to N_{C \to p}(x_1 + \cdots + x_n) \to N_C[x_1 + \cdots + x_n \to p] \to N_{\bar{C}}(p) \to 0.
\end{equation}
This setup is illustrated in Figure~\ref{fig:projection}.

\begin{figure}[h!]
\begin{tikzpicture}
\draw[densely dotted] (1, 4.77) -- (3, -1.6);
\draw[densely dotted] (1, 4.77) -- (-1.37, -2);
\draw[densely dotted] (1, 4.77) -- (2.59, -2.75);
\draw[densely dotted] (1, 4.77) -- (0, -1.31);
\draw[densely dotted] (1, 4.77) -- (1.5, -1.155);
\draw[ultra thick, white] (1.24, 1.975) -- (1.168, 2.8135);
\draw[ultra thick, white] (0.515, 1.825) -- (0.6605, 2.7085);
\draw (2, 5) .. controls (0, 5) and (-1, 2.5) .. (0, 1.5);
\draw (0, 1.5) .. controls (0.5, 1) and (2, 0.75) .. (2, 1.5);
\draw (0.05, 1.55) .. controls (0.5, 2) and (2, 2.25) .. (2, 1.5);
\draw (-0.05, 1.45) .. controls (-1, 0.5) and (1, 0) .. (2, 0);
\filldraw (1, 4.77) circle[radius=0.05];
\draw (-1.5, -1) -- (4.5, -1) -- (3.5, -3) -- (-2.5, -3) -- (-1.5, -1);
\draw (3, -1.6) .. controls (2.75, -0.80375) and (-1.07375, -1.15375) .. (-1.37, -2);
\draw (-1.37, -2) .. controls (-1.66625, -2.84625) and (1, -2.75) .. (2.59, -2.75);
\draw (3, -1.6) .. controls (3.25, -2.39625) and (-0.45, -1.85) .. (-0.85, -1.65);
\filldraw (0.515, 1.825) circle[radius=0.05];
\filldraw (1.24, 1.975) circle[radius=0.05];
\draw[->] (1.24, 1.975) -- (1.168, 2.8135);
\draw[->] (0.515, 1.825) -- (0.6605, 2.7085);
\draw (0.4, 1.96) node{{\tiny \(x_1\)}};
\draw (1.42, 2.08) node{{\tiny \(x_n\)}};
\draw (0.91, 2.02) node{{\tiny \(\cdots\)}};
\draw (1.18, 4.66) node{{\tiny \(p\)}};
\end{tikzpicture}
\caption{Modifications pointing towards \(p\) and projection map}
\label{fig:projection}
\end{figure}
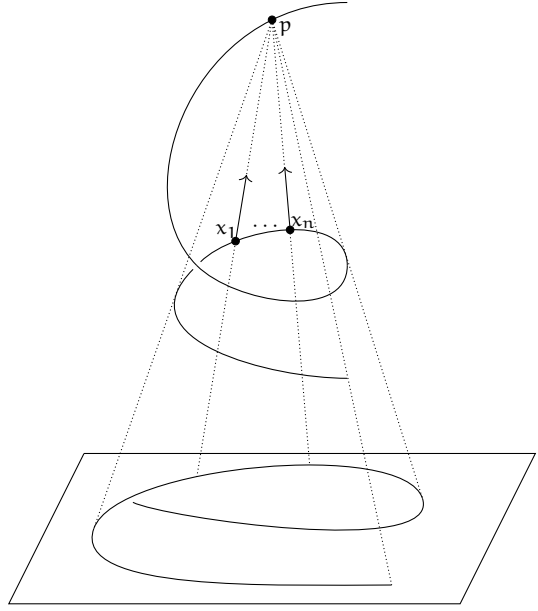

If we get the correct number \(n\) of modifications pointing towards \(p\),
so that the sub and quotient in \eqref{proj} have sufficiently close slope,
then interpolation for \(N_C[x_1 + \cdots + x_n \to p]\)
follows from interpolation for the quotient \(N_{\bar{C}}\).

The main defect of this type of inductive strategy is that it uses up modifications like fuel:
It can only be applied once we already have enough modifications,
and it uses up modifications at each step without creating more.
Moreover, especially when there are only a limited number of different types of modifications,
it can be difficult to get the right number of modifications pointing towards \(p\),
much in the way that you cannot pay \(\$3\) with exact change using only \(\$2\) and \(\$5\) bills.

\bpoint{Miracle: The ``waste product'' of \S \ref{ss:break} can be used as the ``fuel'' required for \S \ref{ss:proj}}
Each time we apply inductive strategies along the lines of \S \ref{ss:break},
we obtain additional modifications.
If we could control these modifications precisely --- so that we can choose how many point towards \(p\) ---
then we could use these modifications to fuel inductive strategies along the lines of \S \ref{ss:proj},
burning up the waste products in the process.

\epoint{Example: Rational space curves of odd degree\label{sss:rsc}}
To illustrate this process in the simplest case,
we prove Theorem~\ref{t:main_normal} for rational curves of odd degree \(d\) in \(\proj^3\).
Peeling off \((d - 3)/2\) one-secant lines as in \S\ref{ss:break}
reduces interpolation for \(N_C\) to interpolation for
\[N_X[2u_1 \to v_1] [2 u_2 \to v_2] \cdots [2 u_{(d - 3)/2} \to v_{(d - 3)/2}],\]
where \(X \subset \proj^3\) is a rational curve of degree \(d - (d - 3)/2 = (d + 3)/2\).
We then specialize \(v_1, v_2, \ldots, v_{(d - 3)/2}\) to a common point \(p \in X\),
so that the resulting bundle fits into an exact sequence
\[0 \to N_{X \to p}(2u_1 + \cdots + 2u_{(d - 3)/2}) \to N_X[2u_1 + \cdots + 2u_{(d - 3)/2} \to p] \to N_{\bar{X}}(p) \to 0,\]
where \(\bar{X}\) is a rational plane curve of degree \((d + 1)/2\)
(and where by \(N_{\bar{X}}\) we mean the normal sheaf of the projection \(\proj^1 = X \to \proj^2\)
rather than the normal bundle of the image).
Both the sub and quotient here are line bundles,
and a straight-forward degree computation shows that both are \(\mathcal{O}_{\proj^1}((3d + 1)/2)\).
Hence, \(N_X[2u_1 + \cdots + 2u_{(d - 3)/2} \to p] \simeq \mathcal{O}_{\proj^1}((3d + 1)/2)^{\oplus 2}\)
satisfies interpolation.

Before discussing the general situation,
we remark that this argument as stated only works when \(d\) is odd.
The case where \(d\) is even is \emph{necessarily} more subtle,
because \(N_C\) does not satisfy interpolation in this case when the characteristic is \(2\).

\epoint{The argument in general}
Despite the rosy picture painted by the above example,
a naive dimension count suggests that the situation is hopeless, due to the third degeneration in \S \ref{ss:break}
(which did not appear in \S\ref{sss:rsc}):
There is only a \((2r - 2)\)-dimensional family of \((r + 1)\)-secant rational curves of degree \(r - 1\),
and it costs \(r - 2\) conditions at each of the \(r + 1\) points for the modification to \(N_X\) to point in a specified direction.
Thus it would be surprising --- in other words, an unlikely intersection on the moduli space
of \((r + 1)\)-secant rational curves of degree \(r - 1\) ---
if we could get more than \((2r - 2)/(r - 2) \approx 2\) of these \(r + 1\) modifications to point towards a specific point \(p\).
The miracle is that this unlikely intersection in fact happens:
As two of the \(r + 1\) points of secancy collide with \(p\),
the curve \(Y\) admits a plethora of degenerations, with \emph{many} different choices for the number of resulting modifications
that point towards the point \(p\).

Of course, showing that appropriate choices can always be made
in the above arguments requires a fairly delicate combinatorial analysis.
Indeed, the presence of the exceptions places a lower bound on the possible complexity of such an argument.
But the realization that this unlikely intersection happens
gives one the faith needed to take the plunge!

\bpoint{Taking stock}  This may seem a baroque induction, but we hope the overview of the counterexamples and the geometry arising from degenerations makes clear how these ideas can lead to a strategy solving the interpolation problem in full generality.  But it is time to take a step back from the subtle details of the proof and review the overall story.

  \section{Summary}

We are led to the study of curves (including Riemann surfaces), and curves in space, from many directions, from Euclid's geometry, to error-correcting codes, to the space of stable maps in Gromov--Witten theory.  When we study curves in space, much as when we study many kinds of mathematical objects, we are led to study moduli spaces of such objects, and compactifications of such spaces correspond to allowing certain kinds of ``degenerate limits" of the objects under consideration.  So to understand smooth curves we are led to instead consider suitable singular curves, where the nature of the "data of non-smoothness" is chosen as necessary to help us prove our theorem.  

The counterexamples to the interpolation problem come because the curves carry with them auxiliary objects (ambient surfaces) which constrain them, so to prove the theorem, we must dodge the counterexamples.    Even though the desired result may be over the complex numbers, strange behavior in characteristic 2 helps guide our search for a proof.  

We hope this problem has given you an excuse for a tour through a number of grand themes in classical and modern geometry which resonate far beyond this particular result.


\begin{thebibliography}{[CLMTB]}

\bibitem[ALY]{p6} A.\ Atanasov, E.\ Larson, and D.\ Yang, {\em Interpolation for normal bundles of general curves},
Mem.\ Amer.\ Math.\ Soc.\ 257 (2019), no.\ 1234, v+105.

\bibitem[B]{4}  G.D. Birkhoff, {\em General mean value value and
    remainder theorems with applications to mechanical differentiation
    and quadrature}, Trans.\ Amer.\ Math.\ Soc.\ {\bf 7} (1906), no.\
  1, 107--136.


\bibitem[CLMTB]{E7}
 A. Castorena, A L\'{o}pez Mart\'{i}n, and M. Teixidor i Bigas, {\em Petri map for vector bundles near good bundles}, J. Pure Appl.\ Algebra {\bf 222} (2018), no.\ 7, 1692--1703.

\bibitem[Ca]{5} A.-L. Cauchy, {\em Cour d'analyse de l'\'Ecole Royale
    Polytecnique; I.re Partie.  Analyse alg\'ebrique}, Chez Debure
  Fr\`eres, Libraires du Roi et de la Biblioth\`eque du Roi, Paris, 1821.



\bibitem[CPJ1]{E15}  K. Cook-Powell and D. Jensen, {\em Components of Brill-Noether loci for curves with fixed gonality}, Michigan Math.\  J. {\bf 71} (2022), no.\ 1, 19--45.

\bibitem[CPJ2]{E16}
 K. Cook-Powell and D. Jensen, {\em Tropical methods in Hurwitz--Brill--Noether theory}, Adv. Math. {\bf 398} (2022), Paper No. 108199.

  \bibitem[Cr]{8} G. Cramer, {\em Introduction \`a l'analyse des lignes
    courbes alg\'ebriques}, Fr\`eres Cramer \& Cl.\ Philbert, Geneva, 1750.

    \bibitem[Eu]{euclid}  Euclid, 
{\em The Elements: Books I--XIII --- Complete and Unabridged},
Sir Thomas Heath trans., Barnes \& Noble, 2006.

\bibitem[EH1]{E10}  D. Eisenbud and J. Harris, {\em Divisors on general curves and cuspidal rational curves}, Invent.\ Math.\ (1983), 371--418.

\bibitem[EH2]{E11}
 D. Eisenbud and J. Harris, {\em Limit Linear Series: Basic Theory}, Invent.\ Math.\ (1986), 337--371.


\bibitem[EH3]{E4}
 D. Eisenbud and J. Harris,  {\em Irreducibility and monodromy of some families of linear series}, Ann.\ Sci.\ Ecole Norm.\ Sup. (4) {\bf 20} (1987), no.\ 1, 65--87.

\bibitem[FL]{E2} W. Fulton and R. Lazarsfeld,  {\em On the connectedness of degeneracy loci and special divisors}, Acta Math. 146 (1981), no. 3-4, 271--283.

\bibitem[G]{E3}  D. Gieseker,  {\em Stable curves and special divisors: Petri's conjecture}, Invent.\ Math.\ {\bf  66} (1982), no.\ 2, 251--275.


\bibitem[GH]{E1} P. Griffiths and J. Harris, {\em On the variety of special linear systems on a general algebraic curve}, Duke Math.\ J.\ {\bf  47} (1980), no.\ 1, 233--272.

 
\bibitem[HH]{hh} R.\ Hartshorne and A.\ Hirschowitz, {\em Smoothing algebraic space curves}, Algebraic geometry, Sitges (Barcelona), 1983, Lecture Notes in Math., vol.\ 1124.

\bibitem[HB]{15}  M. Ch. Hermite and M. Borchardt, {\em Sur la formule
    d'interpolation de Lagrange}, J. Reine Angew.\ Math.\ {\bf 84}
  (1878), 70--79.

\bibitem[JP]{E6}  D. Jensen and S. Payne, {\em Tropical independence I: Shapes of divisors and a proof of the Gieseker--Petri theorem}, Algebra Number Theory {\bf  8} (2014), no.\ 9, 2043--2066.


\bibitem[JR]{JR}
D.\ Jensen and D.\ Ranganathan, {\em Brill--Noether theory for curves of a fixed gonality}, Forum of Math Pi, {\bf 9}, 2021, e1.


\bibitem[K]{E13} G. Kempf,  {\em Schubert methods with an application to algebraic curves}, Pub.\ Math.\ Centrum, Amsterdam (1971).

\bibitem[KL]{E14}
 S. Kleiman and D. Laksov, {\em On the existence of special divisors}, Amer.\ J.\ Math.\ {\bf  94} (1972), 431--436.



    \bibitem[Lag]{17} J. Lagrange, {\em Le\c{c}ons \'el\'ementaires sur les math\'ematiques
        donn\'ees \`a l'\'Ecole Normale en 1795}, in {\em Oeuvres de
        Lagrange}, J. A. Serret, ed., tome 7, 183--288,  Paris, 1867--1892.


\bibitem[Lar]{E17}
 H. Larson, {\em A refined Brill-Noether theory over Hurwitz spaces}, Invent.\ Math.\ {\bf 224} (2021), no.\ 3, 767--790.

\bibitem[LLV]{E9} E. Larson, H. Larson, and I. Vogt, {\em Global Brill--Noether Theory over the Hurwitz Space},  Geom.\ Topol., to appear.

        \bibitem[LV]{lv}  E. Larson and I. Vogt, {\em Interpolation for
          Brill--Noether curves}, Forum of Mathematics Pi, {\bf 11},
        2023, e25.


\bibitem[Laz]{E12}
 R. Lazarsfeld,  {\em Brill--Noether--Petri Without Degenerations}, J.\ Differ.\ Geom. (1986), 299--307.

   \bibitem[M]{mumford} D. Mumford, {\em Further pathologies in algebraic geometry}, Amer.\ J. Math.\ {\bf 84} (1962), 642--648.


\bibitem[O1]{E5}
B. Osserman,  {\em A simple characteristic-free proof of the Brill-Noether theorem}, Bull.\ Braz.\ Math.\ Soc.\ (N.S.)
45 (2014), no.\ 4, 807--818.


\bibitem[O2]{E8}
 B. Osserman, {\em  Connectedness of Brill-Noether loci via degenerations}, Int.\ Math.\ Res.\ Not.\ (2019), no.\ 19, 6162--6178.

\bibitem[Pe]{p5} D.\ Perrin, {\em Courbes passant par \(m\) points g\'en\'eraux de \(\proj^3\)},
M\'em.\ Soc.\ Math.\ France (N.S.) (1987), no.\ 28--29, 138.

\bibitem[Pf]{P} N.\ Pflueger, {\em Brill--Noether varieties of \(k\)-gonal curves}, Adv.\ Math.\ {\bf 312} (2017), 46--63.

\bibitem[R]{p2} Z.\ Ran, {\em Normal bundles of rational curves in projective spaces}, Asian J.\ Math.\ 11 (2007), no.\ 4, 567--608.

   \bibitem[RS]{reedsolomon}  I. S. Reed and G. Solomon, {\em Polynomial codes over certain finite fields}, J. SIAM, Jun.\ 1960 {\bf 8}, no.\ 2, 300--304.

\bibitem[Sa]{p1} G.\ Sacchiero, {\em Normal bundles of rational curves in projective space}, Ann.\ Univ.\ Ferrara Sez.\ VII (N.S.) 26 (1980), 33--40.

\bibitem[St1]{p3} J.\ Stevens, {\em On the number of points determining a canonical curve},
Nederl.\ Akad.\ Wetensch.\ Indag.\ Math.\ 51 (1989), no.\ 4, 485--494.

\bibitem[St2]{p4} J.\ Stevens, {\em On the computation of versal deformations}, vol.\ 82, 1996, Topology, 3, pp.\ 3713--3720. 

\bibitem[V]{murphy} R. Vakil, {\em Murphy's Law in algebraic geometry:  Badly-behaved deformation spaces}, Invent.\  Math.\ {\bf  164} (2006), 569--590. 
   
  \bibitem[W]{27} E. Waring, {\em Problems concerning interpolations},
    Philosophical Transactions of the Royal Society {\bf 69} (1779),
    59--67.

\end{thebibliography}
\end{document}